\documentclass[preprint,11pt]{imsart}

\usepackage{amsthm,amsmath,latexsym,graphicx,verbatim,amsfonts,subfigure}

\startlocaldefs

\newtheorem{theorem}{Theorem}
\newtheorem{lemma}[theorem]{Lemma}
\newtheorem{corollary}[theorem]{Corollary}
\newtheorem{proposition}[theorem]{Proposition}

\def\R{\mathbb{R}}

\def\CC{\mathcal{C}}
\def\HH{\mathcal{H}}
\def\KK{\mathcal{K}}
\def\LL{\mathcal{L}}
\def\NN{\mathcal{N}}
\def\SS{\mathcal{S}}
\def\TT{\mathcal{T}}

\def\alp{\alpha}
\def\del{\delta}
\def\eps{\epsilon}
\def\gam{\gamma}
\def\Gam{\Gamma}
\def\sig{\sigma}
\def\th{\theta}
\def\Th{\Theta}

\def\Ex{\mathbb{E}}
\def\Pr{\mathbb{P}}

\def\argmin{\mathop{\rm argmin}}

\def\bea{\begin{eqnarray*}}
\def\eea{\end{eqnarray*}}

\endlocaldefs

\begin{document}
\addtolength{\baselineskip}{+0.4\baselineskip}

\begin{frontmatter}

% "Title of the paper"
\title{{\sc Technical Report 73}\\
{\sc IMSV, University of Bern}\\[1ex]
Adaptive Confidence Sets for the Optimal Approximating Model}
\runtitle{Confidence Sets for the Best Approximating Model}

% indicate corresponding author with \corref{}
% \author{\fnms{John} \snm{Smith}\corref{}\ead[label=e1]{smith@foo.com}\thanksref{t1}}
\thankstext{t1}{This work was supported by the Swiss National Science Foundation.} 
% \address{line 1\\ line 2\\ printead{e1}}
%\affiliation{Universit\"at Heidelberg}

\author{\fnms{Angelika} \snm{Rohde}\ead[label=e1]{angelika.rohde@math.uni-hamburg.de}}
\address{Universit\"at Hamburg\\ 
Department Mathematik\\
Bundesstra{\ss}e 55\\
D-20146 Hamburg\\
Germany\\
\printead{e1}}
\and
\author{\fnms{Lutz} \snm{D\"umbgen}\ead[label=e2]{duembgen@stat.unibe.ch}}
\address{Universit\"at Bern\\
Institut f\"ur
Mathematische Statistik\\
und Versicherungslehre\\
Sidlerstrasse 5\\
CH-3012 Bern\\
Switzerland\\
\printead{e2}}
\affiliation{Universit\"at Hamburg and Universit\"at Bern}

\runauthor{A. Rohde and L. D\"umbgen}

\begin{abstract}
In the setting of high-dimensional linear models with Gaussian noise, we investigate the possibility of confidence statements connected to model selection. Although there exist numerous procedures for adaptive (point) estimation, the construction of adaptive confidence regions is severely limited (cf.\ Li, 1989). The present paper sheds new light on this gap. We develop exact and adaptive confidence sets for the best approximating model in terms of risk. One of our constructions is based on a multiscale procedure and a particular coupling argument. Utilizing exponential inequalities for noncentral $\chi^2$-distributions, we show that the risk and quadratic loss of all models within our confidence region are uniformly bounded by the minimal risk times a factor close to one.
\end{abstract}



\end{frontmatter}

\newpage

%========================
\section{Introduction}
\label{sec: introduction}
%========================

When dealing with a high dimensional observation vector, the natural question
arises whether the data generating process can be approximated by a model
of substantially lower dimension. Rather than on the true model, the focus is
here on smaller ones which still contain the essential information and allow
for interpretation. Typically, the models under consideration
are characterized by the non-zero components of some parameter
vector. Estimating the true model requires the rather idealistic situation
that each component is either equals zero or has sufficiently modulus: A tiny
perturbation of the parameter vector may result in the biggest model, so
the question about the true model does not seem to be adequate in
general. Alternatively, the model which is optimal in terms of risk appears as a target of many model selection strategies. 
Within a specified class of competing models, this paper is concerned with
confidence regions for those approximating models which are optimal in terms of risk.

Suppose that we observe a random vector $X_n = (X_{in})_{i=1}^n$ with distribution $\NN_n(\th_n,\sig^2 I_n)$ together with an estimator $\hat{\sig}_n$ for the standard deviation $\sig > 0$. Often the signal $\th_n$ represents coefficients of an unknown smooth function with respect to a given orthonormal basis of functions. 

There is a vast amount of literature on point estimation of $\th_n$. For a given estimator $\hat{\th}_n = \hat{\th}_n(X_n,\hat{\sig}_n)$ for $\th_n$, let
$$
	L(\hat{\th}_n, \th_n) \ := \ \|\hat{\th}_n - \th_n\|^2
	\quad\mbox{and}\quad
	R(\hat{\th}_n, \th_n) \ := \ \Ex L(\hat{\th}_n,\th_n)
$$
be its quadratic loss and the corresponding risk, respectively. Here $\|\cdot\|$ denotes the standard Euclidean norm of vectors. Various adaptivity results are known for this setting, often in terms of oracle inequalities. A typical result reads as follows: Let $(\check{\theta}_n^{(c)})_{c \in \CC_n}$ be a family of candidate estimators $\check{\th}_n^{(c)} = \check{\th}_n^{(c)}(X_n)$ for $\th_n$, where $\sig > 0$ is temporarily assumed to be known. Then there exist estimators $\hat{\th}_n$ and constants $A_n, B_n = O(\log(n)^\gam)$ with $\gam \ge 0$ such that for arbitrary $\th_n$ in a certain set $\Th_n \subset \R^n$,
$$
	R(\hat{\th}_n,\th_n)
	\ \le \ A_n \inf_{c \in \CC_n} R(\check{\th}_n^{(c)}, \th_n^{}) + B_n \sigma^2 .
$$
Results of this type are provided, for instance, by Polyak and Tsybakov (1991) and Donoho and Johnstone (1994, 1995, 1998), in the framework of Gaussian model selection by Birg$\acute{\textrm{e}}$ and Massart (2001). The latter article copes in particular with the fact that a model is not necessarily true. Further results of this type, partly in different settings, have been provided by Stone (1984), Lepski et al.\ (1997), Efromovich (1998), Cai (1999, 2002), to mention just a few.

By way of contrast, when aiming at adaptive confidence sets one faces severe limitations. Here is a result of Li (1989), slightly rephrased: Suppose that $\Th_n$ contains a closed Euclidean ball $B(\th_n^o, c n^{1/4})$ around some vector $\theta_n^o \in \R^n$ with radius $c n^{1/4} > 0$. Still assuming $\sig$ to be known, let $\hat{D}_n = \hat{D}_n(X_n) \subset \Th_n$ be a $(1 - \alpha)$-confidence set for $\th_n \in \Th_n$. Such a confidence set may be used as a test of the (Bayesian) null hypothesis that $\th_n$ is uniformly distributed on the sphere $\partial B(\th_n^o, c n^{1/4})$ versus the alternative that $\th_n = \th_n^o$: We reject this null hypothesis at level $\alpha$ if $\|\eta - \th_n^o\| < c n^{1/4}$ for all $\eta \in \hat{D}_n$. Since this test cannot have larger power than the corresponding Neyman-Pearson test,
\bea
	\Pr_{\th_n^o}^{} \biggl( \sup_{\eta \in \hat{D}_n} \|\eta - \th_n^o\|_n < c n^{1/4} \biggr)
	& \le & \Pr \Bigl( S_n^2 \le \chi_{n; \alpha}^2(c^2 n^{1/2} / \sig^2) \Bigr)
		\quad	(S_n^2 \sim \chi_n^2) \\
	& = & \Phi \Bigl( \Phi^{-1}(\alpha) + 2^{-1/2} c^2 / \sig^2 \Bigr) + o(1) ,
\eea
where $\chi_{n;\alp}^2 (\del^2)$ stands for the $\alp$-quantile of the noncentral chi-squared distribution with $n$ degrees of freedom and noncentrality parameter $\del^2$. Throughout this paper, asymptotic statements refer to $n \to \infty$. The previous inequality entails that no reasonable confidence set has a diameter of order $o_p(n^{1/4})$ uniformly over the parameter space $\Th_n$, as long as the latter is sufficiently large. Despite these limitations, there is some literature on confidence sets in the present or similar settings; see for instance Beran (1996, 2000), Beran and D\"{u}mbgen (1998) and Genovese and Wassermann (2005).

Improving the rate of $O_p(n^{1/4})$ is only possible via additional constraints on $\th_n$, i.e.\ considering substantially smaller sets $\Th_n$. For instance, Baraud (2004) developed nonasymptotic confidence regions which perform well on finitely many linear subspaces. Robins and van der Vaart (2006) construct confidence balls via sample splitting which adapt to some extent to the unknown ``smoothness'' of $\th_n$. In their context, $\Th_n$ corresponds to a Sobolev smoothness class with given parameter $(\beta, L)$. However, adaptation in this context is possible only within a range $[\beta, 2 \beta]$. Independently, Cai and Low (2006) treat the same problem in the special case of the Gaussian white noise model, obtaining the same kind of adaptivity in the broader scale of Besov bodies. Other possible constraints on $\th_n$ are so-called shape constraints; see for instance Cai and Low (2007), D\"umbgen (2003) or Hengartner and Stark (1995).

The question is whether one can bridge this gap between confidence sets and point estimators. More precisely, we would like to understand the possibility of adaptation for point estimators in terms of some confidence region for the set of all optimal candidate estimators $\check{\th}_n^{(c)}$. That means, we want to construct a confidence region $\hat{\KK}_{n,\alpha} = \hat{\KK}_{n,\alpha}(X_n,\hat{\sig}_n) \subset \CC_n$ for the set
\begin{align*}
	\KK_n(\th_n) \, :=& \, \argmin_{c\in\CC_n}R(\check{\th}_n^{(c)}) \\
	=& \, \Big\{c\in\CC_n: \, R(\check{\th}_n^{(c)}, \th_n)\leq R(\check{\th}_n^{(c')}, \th_n)\, \text{for all}\, c'\in\CC_n\Big\}
\end{align*}
such that for arbitrary $\th_n \in \R^n$,
\begin{equation}
	\Pr_{\theta_n} \Bigl( \KK_{n}(\th_n) \subset \hat{\KK}_{n,\alpha} \Bigr) \ \ge \ 1 - \alp
	\label{equ: confidence level of KKhat}
\end{equation}
and
\begin{equation}
	\left.\begin{array}{c}
		\displaystyle
		\max_{c \in \hat{\KK}_{n,\alpha}} \, R(\check{\th}_n^{(c)}, \th_n^{}) \\
		\displaystyle
		\max_{c \in \hat{\KK}_{n,\alpha}} \, L(\check{\th}_n^{(c)}, \th_n^{})
	\end{array}\right\}
	\ = \ O_p(A_n) \min_{c \in \CC_n} \, R(\check{\th}_n^{(c)}, \th_n^{}) + O_p(B_n) \sig^2 .
	\label{equ: quality of KKhat}
\end{equation}
Solving this problem means that statistical inference about differences in the performance of estimators is possible, although inference about their risk and loss is severely limited. In some settings, selecting estimators out of a class of competing estimators entails estimating implicitly an unknown regularity or smoothness class for the underlying signal $\th_n$. Computing a confidence region for good estimators is particularly suitable in situations in which several good candidate estimators fit the data equally well although they look different. This aspect of exploring various candidate estimators is not covered by the usual theory of point estimation.

Note that our confidence region $\hat{\KK}_{n,\alpha}$ is required to contain the whole set $\KK_n(\th_n)$, not just one element of it, with probability at least $1 - \alp$. The same requirement is used by Futschik (1999) for inference about the argmax of a regression function.

The remainder of this paper is organized as follows. For the reader's convenience our approach is first described in a simple toy model in Section~\ref{sec: toy model}. In Section~\ref{sec: nested} we develop and analyze an explicit confidence region $\hat{\KK}_{n,\alpha}$ related to $\CC_n := \{0,1,\ldots,n\}$ with candidate estimators
$$
	\check{\th}_n^{(k)} \ := \ \bigl( 1\{i \le k\} X_{in} \bigr)_{i=1}^n .
$$
These correspond to a standard nested sequence of approximating models. Section~\ref{sec: general candidates} discusses richer families of candidate estimators.

All proofs and auxiliary results are deferred to Sections~\ref{sec: proofs} and \ref{sec: appendix}.

%======================
\section{A toy problem}
\label{sec: toy model}
%======================

Suppose we observe a stochastic process $Y = (Y(t))_{t \in [0,1]}$, where
$$
	Y(t) \ = \ F(t) + W(t) ,	\quad	t \in [0,1] ,
$$
with an unknown fixed continuous function $F$ on $[0,1]$ and a Brownian motion $W = (W(t))_{t\in [0,1]}$. We are interested in the set
$$
	\SS(F) \ := \ \argmin_{t \in [0,1]} F(t) .
$$
Precisely, we want to construct a $(1 - \alpha)$-confidence region $\hat{\SS}_{\alpha} = \hat{\SS}_{\alpha}(Y) \subset [0,1]$ for $\SS(F)$ in the sense that
\begin{equation}
	P \bigl( \SS(F) \subset \hat{\SS}_{\alpha}\bigr) \ \ge \ 1 - \alpha ,
	\label{Level0}
\end{equation}
regardless of $F$. To construct such a confidence set we regard $Y(s) - Y(t)$ for arbitrary different $s, t \in [0,1]$ as a test statistic for the null hypothesis that $s \in \SS(F)$, i.e.\ large values of $Y(s) - Y(t)$ give evidence for $s \not\in \SS(F)$.

A first naive proposal is the set
$$
	\hat{\SS}_{\alpha}^{\rm naive} \ := \ \Bigl\{ s \in [0,1] : Y(s) \le
        \min_{[0,1]} Y + \kappa_\alpha^{\rm naive} \Bigr\}
$$
with $\kappa_\alpha^{\rm naive}$ denoting the $(1 - \alp)$-quantile of $\max_{[0,1]} W - \min_{[0,1]} W$.

Here is a refined version based on results of D\"umbgen and Spokoiny (2001): Let $\kappa_\alpha$ be the $(1 - \alpha)$-quantile of
\begin{equation}
	\label{eqn: toy test statistic}
	\sup_{s,t \in [0,1]} \biggl( \frac{|W(s) - W(t)|}{\sqrt{|s-t|}} - \sqrt{2 \log(e/\arrowvert s-t\arrowvert)} \biggr).
\end{equation}
Then constraint (\ref{Level0}) is satisfied by the confidence region
$\hat{\SS}_{\alpha}$ which consists of all $s \in [0,1]$ such that
$$
	Y(s) \le Y(t) + \sqrt{|s-t|} \Bigl(\sqrt{2 \log(e/\arrowvert s-t\arrowvert)}
		\,+\, \kappa_\alpha \Bigr)
		\ \text{for all} \ t \in [0,1] .
$$

To illustrate the power of this method, consider for instance a sequence of functions $F = F_n = c_n F_o$ with positive constants $c_n \to \infty$ and a fixed continuous function $F_o$ with unique minimizer $s_o$. Suppose that
$$
	\lim_{t \to s_o} \frac{F_o(t) - F_o(s_o)}{|t - s_o|^\gamma} \ = \ 1
$$
for some $\gamma > 1/2$. Then the naive confidence region satisfies only
\begin{equation}
\label{eq: performance naive S}
	\max_{t \in \hat{\SS}_{\alpha}^{\rm naive}} \, |t - s_o|
	\ = \ O_p \bigl( c_n^{-1/\gam} \bigr) ,
\end{equation}
whereas
\begin{equation}
\label{eq: performance advanced S}
	\max_{t \in \hat{\SS}_{\alpha}} \, |t - s_o|
	\ = \ O_p \Bigl( \log(c_n)^{1/(2\gam - 1)} c_n^{-2/(2\gam - 1)} \Bigr) .
\end{equation}

%===========================================================
\section{Confidence regions for nested approximating models}
\label{sec: nested}
%===========================================================

As in the introduction let $X_n = \th_n + \eps_n$ denote the $n$-dimensional observation vector with $\theta_n\in\R^n$ and $\eps_n \sim \mathcal{N}_n(0,\sigma^2 I_n)$. For any candidate estimator $\check{\th}_n^{(k)} = \bigl( 1\{i \le k\} X_{in} \bigr)_{i=1}^n$ the loss is given by 
$$
	L_n(k) := L(\check{\th}_n^{(k)},\th_n)
	\ = \ \sum_{i=k+1}^n \th_{in}^2 + \sum_{i=1}^k (X_{in} - \th_{in})^2
$$
with corresponding risk
$$
	R_n(k) := R(\check{\th}_n^{(k)},\th_n)
	\ = \ \sum_{i=k+1}^n \th_{in}^2 + k \sig^2 .
$$
Model selection usually aims at estimating a candidate estimator which is optimal in terms of risk. Since the risk depends on the unknown signal and therefore is not available, the selection procedure minimizes an unbiased risk estimator instead. In the sequel, the bias-corrected risk estimator for the candidate $\check{\th}_n^{(k)}$ is defined as
$$
	\hat{R}_n(k) \ := \ \sum_{i=k+1}^n (X_{in}^2 - \hat{\sig}_n^2) +k \hat{\sig}_n^2 ,
$$
where $\hat{\sigma}_n^2$ is a variance estimator satisfying the subsequent condition.

\medskip
\noindent
\textbf{(A)} $\hat{\sig}_n^2$ and $X_n$ are stochastically independent with
$$
	\frac{m \hat{\sig}_n^2}{\sig^2} \ \sim \ \chi_{m}^2 ,
$$
where $1 \le m = m_n \le \infty$ with $m = \infty$ meaning that $\sig$ is
known, i.e. $\hat{\sigma}_n^2\equiv \sigma^2$. For asymptotic statements, it
is generally assumed that 
$$
	\beta_n^2 \ := \ \frac{2n}{m_n} \ = \ O(1)
$$ 
unless stated otherwise.

\paragraph{Example.}
Suppose that we observe $Y = M \eta + \del$ with given design matrix $M \in
\R^{(n+m)\times n}$ of  rank $n$, unknown parameter vector $\eta \in \R^n$ and
unobserved error vector $\del \sim \NN_{n+m}(0, \sig^2 I_{n+m})$. Then the previous
assumptions are satisfied by $X_n := (M^\top M)^{1/2} \hat{\eta}$ with
$\hat{\eta} := (M^\top M)^{-1} M^\top Y$ and $\hat{\sig}_n^2 := \|Y - M
\hat{\eta}\|^2 / m$, where $\th_n := (M^\top M)^{1/2} \eta$.

\medskip

Important for our analysis is the behavior of the centered and rescaled difference process 
$\hat{D}_n = \bigl( \hat{D}_n(j,k) \bigr)_{0\le j < k \le n}$ with
\begin{eqnarray*}
	\hat{D}_n(j,k)
	& := & 
		\frac{\hat{R}_n(j) - \hat{R}_n(k) - R_n(j) + R_n(k)}{\hat{\sigma}_n^2\big(4\Arrowvert\theta_n/\sigma\Arrowvert^2+2n\big)^{1/2}} \\
	& = & \frac{\sum_{i=j+1}^k (X_{in}^2 - \sigma^2 - \theta_{in}^2)
			- 2 (k - j) (\hat{\sigma}^2 - \sigma^2)}{\hat{\sigma}_n^2\big(4\Arrowvert\theta_n/\sigma\Arrowvert^2+2n\big)^{1/2}} .
\end{eqnarray*}
One may also write $\hat{D}_n(j,k) = (\hat{\sigma}_n/\sigma)^{-2} \bigl( D_n(j,k) + V_n(j,k) \bigr)$ with
\begin{align}
\label{eq: Definition D_n}
	D_n(j,k)\, & :=\,  \frac{1}{\sqrt{4\Arrowvert\theta_n/\sigma\Arrowvert^2+2n}}\sum_{i=j+1}^k
		\Bigl( 2 (\theta_{in}/\sigma) (\eps_{in}/\sigma) + (\eps_{in}/\sigma)^2 - 1 \Bigr) , \\
\label{eq: Definition V_n}
	V_n(j,k)\, & := \, \big(4\Arrowvert\theta_n/\sigma\Arrowvert^2+2n\big)^{-1/2}2 (k - j) (1 - \hat{\sigma}^2/\sigma^2)
\end{align}
This representation shows that the distribution of $\hat{D}_n$ depends on the degrees of freedom, $m$, and the unknown ``signal-to-noise vector" $\theta_n/\sigma$. The process $D_n$ consists of partial sums of the independent, but in general non-identically distributed random variables $2 (\theta_{in}/\sigma) (\eps_{in}/\sigma) + (\eps_{in}/\sigma)^2 - 1$. The standard deviation of $D_n(j,k)$ is given by
$$
	\tau_n(j,k)
	\ := \ \frac{1}{\sqrt{4\Arrowvert\theta_n/\sigma\Arrowvert^2+2n}}
		\Bigl( \sum_{i=j+1}^k \bigl( 4 \theta_{in}^2/\sigma^2 + 2 \bigr) \Bigr)^{1/2} .
$$
Note that $\tau_n(0,n)=1$ by construction. 
To imitate the more powerful confidence region of Section~\ref{sec: toy model} based on the multiscale approach, one needs a refined analysis of the increment process $\hat{D}_n$. Since this process does not have subgaussian tails, the standardization is more involved than the correction in (\ref{eqn: toy test statistic}).

\begin{theorem}
\label{thm: approximation 2}
Define $\Gamma_n(j,k) := \big(2 \log\big(e/\tau_n(j,k)^2\big)^{1/2}$ for $0 \le j < k \le n$. Then
$$
	\sup_{0 \le j < k \le n} \frac{|\hat{D}_n(j,k)|}{\tau_n(j,k)}
	\ \le \ \sqrt{32} \, \log n + O_p(1) ,
$$
and for any fixed $c > 2$,
$$
	\hat{d}_n := \sup_{0 \le j < k \le n}
		\Biggl( \frac{|\hat{D}_n(j,k)|}{\tau_n(j,k)}
			\,-\, \Gamma_n(j,k)\, -\,  \frac{c\cdot\Gamma_n(j,k)^2}{\big(4\Arrowvert\theta_n/\sigma\Arrowvert^2+2n\big)^{1/2}\tau_n(j,k)} \Biggr)^+
$$
is bounded in probability. In case of $\|\theta_n\|^2 = O(n)$, $\LL(\hat{d}_n)$ is weakly approximated by the law of
$$
	\delta_n := \sup_{0 \le j < k \le n}
		\biggl( \frac{|\Delta_n(j,k)|}{\tau_n(j,k)}
			- \Gamma_n(j,k) \biggr)^+,
$$
where
$$
	\Delta_n(j,k) 
	\ = \ W(\tau_n(0,k)^2) - W(\tau_n(0,j)^2)
		-  \frac{2\beta_n(k - j)}{\sqrt{n} \, \big(4\Arrowvert\theta_n/\sigma\Arrowvert^2+2n\big)^{1/2}} \, Z
$$
with a standard Brownian motion $W$ and a random variable $Z \sim \mathcal{N}(0,1)$ independent of $W$.
\end{theorem}

The limiting distribution indicates that the additive correction term in the definition of $\hat{d}_n$ cannot be chosen essentially smaller. It will play a crucial role for the efficiency of the confidence region. 

To construct a confidence set for $\KK_n(\theta_n)$ by means of $\hat{d}_n$, we are facing the problem that the auxiliary function $\tau_n(\cdot,\cdot)$ depends on the unknown signal-to-noise vector $\th_n/\sig$. In fact, knowing $\tau_n$ would imply knowledge of $\KK_n(\th_n)$ already. A natural approach is to replace the quantities which are dependent on the unknown parameter by suitable estimates. A common estimator of the variance $\tau_n(j,k)^2$, $j < k$, is given by
$$
	\hat{\tau}_n(j,k)^2
	\ := \ \bigg\{\sum_{i=1}^n
		\Bigl( 4 (X_{in}^2/\hat{\sigma}_n^2 - 1) + 2 \Bigr)\biggr\}^{-1}\sum_{i=j+1}^k
		\Bigl( 4 (X_{in}^2/\hat{\sigma}_n^2 - 1) + 2 \Bigr) .
$$
However, using such an estimator does not seem to work since
$$
	\sup_{0 \le j < k \le n} \, \Bigl| \frac{\hat{\tau}_n(j,k)}{\tau_{n}(j,k)} - 1 \Bigr|
	\ \not\longrightarrow_{p} \ 0
$$
as $n$ goes to infinity. This can be verified by noting that the (rescaled) numerator of  $\bigl(\hat{\tau}_n(j,k)^2\bigr)_{0 \le j < k\le n}$ is, up to centering, essentially of the same structure as the rescaled difference process $\hat{D}_n$ itself.

%-------------------------------------------------------
\subsection*{The least favourable case of constant risk}
%-------------------------------------------------------

The problem of estimating the set $\arg\min_k R_n(k)$ can be cast into our toy model where $Y(t)$, $F(t)$ and $W(t)$ correspond to $\hat{R}_n(k)$, $R_n(k)$ and the difference $\hat{R}_n(k) - R_n(k)$, respectively. One may expect that the more distinctive the global minima are, the easier it is to identify their location. Hence the case of constant risks appears to be least favourable, corresponding to a signal
$$
	\th_n^* \ := \ \bigl( \pm \sig \bigr)_{i=1}^n ,
$$
In this situation, each candidate estimator $\check{\th}_n^{(k)}$ has the same risk of $n\sigma^2$.

A related consideration leading to an explicit procedure is as follows: For fixed indices $0 \le j < k \le n$,
$$
	R_n(j) - R_n(k) \ = \ \sum_{i=j+1}^k \th_{in}^2 - (k - j) \sigma^2 ,
$$
and if Assumption (A) is satisfied, the statistic
$$
	T_{jkn}
	\ := \ \frac{\sum_{i=j+1}^k X_{in}^2}{(k - j) \hat{\sig}_n^2}
	\ = \ 2 - \frac{\hat{R}_n(k) - \hat{R}_n(j)}{(k - j) \hat{\sig}_n^2}
$$
has a noncentral (in the numerator) $F$-distribution
$$
	F_{k-j, m}^{} \biggl( \frac{\sum_{i=j+1}^k \th_{in}^2}{\sig^2} \biggr)
		= F_{k-j,m}^{} \biggl( k - j + \frac{R_n(j) - R_n(k)}{\sig^2} \biggr)
$$
with $k-j$ and $m$ degrees of freedom. Thus large or small values of $T_{jkn}$ give evidence for $R_n(j)$ being larger or smaller, respectively, than $R_n(k)$. Precisely,
$$
	\LL_{\th_n}(T_{jkn}) \ \begin{cases}
		\le_{\rm st.} \ \LL_{\th_n^*}(T_{jkn}) & \mbox{whenever } j \in \KK_n(\th_n) , \\
		\ge_{\rm st.} \ \LL_{\th_n^*}(T_{jkn}) & \mbox{whenever } k \in \KK_n(\th_n) .
	\end{cases}
$$
Note that this stochastic ordering remains valid if $\hat{\sig}_n^2$ is just
independent from $X_n$, i.e. also under the more general requirement of the
remark at the end of this section. Via suitable coupling of Poisson mixtures of central $\chi^2$-distributed random variables, this observation is extended to a coupling for the whole process $\big(T_{jkn}\big)_{0\leq j<k\leq n}$:

\begin{proposition}[Coupling]
\label{prop: coupling}
For any $\th_n \in \R^n$ there exists a probability space with random variables $\bigl( \tilde{T}_{jkn} \bigr)_{0 \le j < k \le n}$ and $\bigl( \tilde{T}^*_{jkn} \bigr)_{0 \le j < k \le n}$ such that
\bea
	\LL \Bigl( \bigl( \tilde{T}_{jkn} \bigr)_{0 \le j < k \le n} \Bigr)
	& = & \LL_{\th_n}^{} \Bigl( \bigl( T_{jkn} \bigr)_{0 \le j < k \le n} \Bigr) , \\
	\LL \Bigl( \bigl( \tilde{T}^*_{jkn} \bigr)_{0 \le j < k \le n} \Bigr)
	& = & \LL_{\th_n^*}^{} \Bigl( \bigl( T_{jkn} \bigr)_{0 \le j < k \le n} \Bigr) ,
\eea
and for arbitrary indices $0 \le j < k \le n$,
$$
	\tilde{T}_{jkn} \ \begin{cases}
		\le \ \tilde{T}_{jkn}^* & \mbox{whenever } j \in \KK_n(\th_n) , \\
		\ge \ \tilde{T}_{jkn}^* & \mbox{whenever } k \in \KK_n(\th_n) .
	\end{cases}
$$
\end{proposition}

As a consequence of Proposition~\ref{prop: coupling}, we can define a confidence set for $\KK_n(\theta_n)$, based on this least favourable case. Let $\kappa_{n,\alpha}$ denote the $(1-\alpha)$-quantile of $\LL_{\theta_n^*}(\hat{d}_n)$, where for simplicity $c:=3$ in the definition of $\hat{d}_n$. Note also that $\tau_n(j,k)^2 = (k - j)/n$ in case of $\theta_n = \theta_n^*$. Motivated by the procedure in Section~\ref{sec: toy model} and Theorem~\ref{thm: approximation 2}, we define  
\begin{eqnarray}
	\hat{\KK}_{n,\alpha}
	& := & \Bigl\{ j : \,
		\hat{R}_n(j) \le \hat{R}_n(k) + \hat{\sigma}_n^2 |k - j| c_{jkn}
		\ \text{for all} \ k \ne j \Bigr\}
	\label{equ: definition KKhat}  \\
	& \, = & \bigl\{ j : \,
		T_{ijn} \ge 2 - c_{ijn} \ \text{for all} \ i < j,
	\nonumber \\
	&& \qquad\qquad\qquad\qquad\quad
		T_{jkn} \le 2 + c_{jkn} \ \text{for all} \ k > j \bigr\}
	\nonumber
\end{eqnarray}
with 
$$
	c_{jkn} = c_{jkn,\alpha}
	\ := \ \sqrt{\frac{6}{|k-j|}} \biggl( \Gam \Bigl(
        \frac{ k-j}{n} \Bigr) + \kappa_{n, \alpha} \bigg)
		+ \frac{3}{|k-j|} \Gam \Bigl(
                \frac{k-j}{n} \Bigr)^2.
$$

\begin{theorem}
\label{thm: oracle nested}
Let $(\theta_n)_{n\in\mathbb{N}}$ be arbitrary. With $\hat{\KK}_{n,\alpha}$ as defined above,
$$
	\Pr_{\theta_n}\Bigl(\KK_n(\theta_n)\not\subset\hat{\KK}_{n,\alpha}\Bigr)
	\ \le \ \alpha.
$$
In case of $\beta_n \to 0$ (i.e.\ $n /m \to 0$), the critical values $\kappa_{n, \alpha}$ converge to the critical value $\kappa_\alp$ introduced in Section~\ref{sec: toy model}. In general, $\kappa_{n, \alpha} = O(1)$, and the confidence regions $\hat{\KK}_{n,\alpha}$ satisfy the oracle inequalities
\begin{align}
	\max_{k \in \hat{\KK}_{n,\alpha}} R_n(k)
	\ \le \ \min_{j \in \CC_n} \, R_n(j) \
		& + \ \bigl( 4\sqrt{3} + o_p(1) \bigr)
			\sqrt{\sigma^2\log(n) \min_{j\in\CC_n}R_n(j)}
	\label{equ: quality of KKhat1} \\
		& + \ O_p\bigl(\sigma^2 \log n\bigr)
	\nonumber
\intertext{and}
	\max_{k \in \hat{\KK}_{n,\alpha}} L_n(k)
	\ \le \ \min_{j \in \CC_n} \, L_n(j) \
		& + \ O_p \Bigl( \sqrt{ \sigma^2 \log(n) \min_{j\in\CC_n}L_n(j) } \Bigr)
	\label{equ: quality of KKhat2} \\
		& + \ O_p\bigl(\sigma^2 \log n\bigr) .
	\nonumber
\end{align}
\end{theorem}

\paragraph{Remark (Dependence on $\alpha$)}
The proof reveals a refined version of the bounds in Theorem \ref{thm: oracle nested} in case of signals $\theta_n$ such that
$$
	\log(n)^3 \ = \ O \Bigl( \min_{j \in \CC_n} \, R_n(j) \Bigr) .
$$
Let $0 < \alpha(n) \to 0$ such that $\kappa_{n,\alpha(n)}^6 = O \Bigl( \min_{j \in \CC_n} \, R_n(j) \Bigr)$. Then
\bea
	\max_{k \in \hat{\KK}_{n,\alpha}} \, R_n(k)
	& \le & \min_{j \in \CC_n} \, R_n(j) \\
	&& + \ \Bigl( 4 \sqrt{3} \sqrt{\log n} + 2 \sqrt{6} \, \kappa_{n,\alpha}
		+ O_p(1) \Bigr) \sqrt{\sigma^2 \min_{j\in\CC_n}R_n(j)}
\eea
uniformly in $\alpha \ge \alpha(n)$.

\paragraph{Remark (Variance estimation)}
Instead of Condition (A), one may require more generally that
$\hat{\sig}_n^2$ and $X_n$ are independent with
$$
	\sqrt{n} \Bigl( \frac{\hat{\sig}_n^2}{\sig^2} - 1 \Bigr) \ \to_D^{} \ \NN(0, \beta^2)
$$
for a given $\beta \ge 0$. This covers, for instance, estimators used in connection with wavelets. There $\sig$ is estimated by the median of some high frequency wavelet coefficients divided by the normal quantile $\Phi^{-1}(3/4)$. Theorem~\ref{thm: approximation 2} continues to hold, and the coupling extends to this situation, too, with $S^2$ in the proof being distributed as $n\hat{\sigma}_n^2$. Under this assumption on the external variance estimator, the confidence region $\hat{\KK}_{n,\alpha}$, defined with $m := \lfloor 2n/\beta^2\rfloor$, is at least asymptotically valid and satisfies the above oracle inequalities as well.

%====================================================================
\section{Confidence sets in case of larger families of candidates}
\label{sec: general candidates}
%====================================================================

The previous result relies strongly on the assumption of nested models. It is
possible to obtain confidence sets for the optimal approximating models in a more general setting, albeit the resulting oracle property is not as strong as in the nested case. In particular, we can no longer rely on a coupling result but need a different construction. For the reader's convenience, we focus on the case of known $\sigma$, i.e.\ $m = \infty$; see also the remark at the end of this section.

Let $\CC_n$ be a family of index sets $C \subset \{1,2,\ldots,n\}$ with candidate estimators
$$
	\check{\th}^{(C)} \ := \ \bigl( 1\{i \in C\} X_{in} \bigr)_{i=1}^n 
$$
and corresponding risks
$$
	R_n(C) := R(\check{\th}^{(C)}, \th_n)
	\ = \ \sum_{i \not\in C} \th_{in}^2 + |C| \sigma^2 ,
$$
where $|S|$ denotes the cardinality of a set $S$. For two index sets $C$ and $D$,
$$
	\sigma^{-2} \bigl( R_n(D) - R_n(C) \bigr)
	\ = \ \delta_n^2(C \setminus D) - \delta_n^2(D \setminus C) + |D| - |C|
$$
with the auxiliary quantities
$$
	\delta_n^2(J) \ := \ \sum_{i \in J} \th_{in}^2/\sig^2 ,	\quad J \subset \{1,2,\ldots,n\} .
$$
Hence we aim at simultaneous $(1 - \alpha)$-confidence intervals for these noncentrality parameters $\delta_n(J)$, where $J \in \mathcal{M}_n := \{D\setminus C : C,D \in \CC_n\}$. To this end we utilize the fact that
$$
	T_n(J) \ := \ \frac{1}{\sigma^2} \sum_{i\in J} X_{in}^2
$$
has a $\chi_{|J|}^2(\delta_n^2(J))$-distribution. We denote the distribution function of $\chi_k^2(\delta^2)$ by $F_k(\cdot \mid \delta^2)$. Now let $M_n := |\mathcal{M}_n| - 1 \le |\CC_n|(|\CC_n| - 1)$, the number of nonvoid index sets $J \in \mathcal{M}_n$. Then with probability at least $1 - \alpha$,
\begin{equation}
	\alpha/(2M_n)
	\le F_{|J|} \bigl( T_n(J) \,\big|\, \delta_n^2(J) \bigr)
	\le 1 - \alpha/(2M_n)
	\quad\text{for} \ \emptyset \ne J \in \mathcal{M}_n .
	\label{CC.0}
\end{equation}
Since $F_{|J|}(T_n(J) \mid \delta^2)$ is strictly decreasing in $\delta^2$ with limit $0$ as $\delta^2 \to \infty$, (\ref{CC.0}) entails the simultaneous $(1 - \alpha)$-confidence intervals $\bigl[ \hat{\delta}_{n,\alpha,l}^2(J), \hat{\delta}_{n,\alpha,u}^2(J) \bigr]$ for all parameters $\delta_n^2(J)$ as follows: We set $\hat{\delta}_{n,\alpha,l}^2(\emptyset) := \hat{\delta}_{n,\alpha,u}^2(\emptyset) := 0$, while for nonvoid $J$,
\begin{eqnarray}
	\hat{\delta}_{n,\alpha,l}^2(J)
	& := & \min \Bigl\{ \delta^2 \ge 0 :
		F_{|J|} \bigl( T_n(J) \,\big|\, \delta^2 \bigr)
		\le 1 - \alpha/(2M_n) \Bigr\} ,
	\label{CC.1L} \\
	\hat{\delta}_{n,\alpha,u}^2(J)
	& := & \max \Bigl\{ \delta^2 \ge 0 :
		F_{|J|} \bigl( T_n(J) \,\big|\, \delta^2 \bigr)
		\ge \alpha/(2M_n) \Bigr\} .
	\label{CC.1U}
\end{eqnarray}
By means of these bounds, we may claim with confidence $1 - \alpha$ that for arbitrary $C, D \in \CC_n$ the normalized difference $(n/\sig^2) \bigl( R_n(D) - R_n(C) \bigr)$ is at most $\hat{\delta}_{n,\alp,u}^2(C \setminus D) - \hat{\delta}_{n,\alp,l}^2(D \setminus C) + |D| - |C|$. Thus a $(1 - \alp)$-confidence set for $\KK_n(\th_n) = \argmin_{C \in \CC_n} \, R_n(C)$ is given by
$$
	\hat{\KK}_{n,\alp}
	\ := \ \Bigl\{ C \in \CC_n :
		\hat{\delta}_{n,\alp,u}^2(C \setminus D) - \hat{\delta}_{n,\alp,l}^2(D \setminus C) + |D| - |C| \ge 0
		\ \text{for all} \ D \in \CC_n \Bigr\} .
$$
These confidence sets $\hat{\KK}_{n,\alp}$ satisfy the following oracle inequalities:

\begin{theorem}
\label{thm: oracle general}
Let $(\theta_n)_{n\in\mathbb{N}}$ be arbitrary, and suppose that $\log |\CC_n| = o(n)$. Then
\begin{align*}
	\max_{C \in \hat{\KK}_{n,\alpha}} \, R_n(C)
	\ \le \ \min_{D \in \CC_n} \, R_n(D) \
		& + \ O_p \biggl( \sqrt{ \sigma^2 \log(\arrowvert\CC_n\arrowvert)
			\min_{D \in \CC_n} \, R_n(D) } \biggr) \\
		& + \ O_p \big( \sigma^2 \log \arrowvert \CC_n\arrowvert\big)
\intertext{and}
	\max_{C \in \hat{\KK}_{n,\alpha}} \, L_n(C)
	\ \le \ \min_{D \in \CC_n} \, L_n(D) \
		& + \ O_p \biggl( \sqrt{ \sigma^2 \log(\arrowvert\CC_n\arrowvert)
			\min_{D \in \CC_n} \, L_n(D) } \biggr) \\
		& + \ O_p \big( \sigma^2\log \arrowvert\CC_n\arrowvert\big). 
\end{align*}
\end{theorem}

\paragraph{Remark.}
The upper bounds in Theorem \ref{thm: oracle general} are of the form
$$
	\rho_n \Bigl( 1
		+ O_p \Bigl( \sqrt{\sigma^2 \log(\arrowvert\CC_n\arrowvert)/\rho_n}
		+ \sigma^2\log(\arrowvert\CC_n\arrowvert)/\rho_n \Bigr) \Bigr)
$$
with $\rho_n$ denoting minimal risk or minimal loss. Thus Theorem~\ref{thm: oracle general} entails that the maximal risk (loss) over $\hat{\KK}_{n,\alpha}$ exceeds the minimal risk (loss) only by a factor close to one, provided that the minimal risk (loss) is substantially larger than $\sigma^2\log\arrowvert\CC_n\arrowvert$.

\paragraph{Remark (Suboptimality in case of nested models)}
In case of nested models, the general construction is suboptimal in the factor of the leading (in most cases) term $\sqrt{\min_j R_n(j)}$. Following the proof carefully and using that $\sigma^2\log\arrowvert\CC_n\arrowvert=2\sigma^2\log n\,+\, O(1)$ in this special setting, one may verify that
\begin{align*}
	\max_{k \in \hat{\KK}_{n,\alpha}} \, R_n(k)
	\ \le \ \min_{j \in \CC_n} \, R_n(j) \
		& + \ \bigl( 4\sqrt{8} + o_p(1) \bigr)
			\sqrt{ \sigma^2 \log(n) \min_{j \in \CC_n} \, R_n(j) } \\
		& + \ O_p \big( \sigma^2 \log n\big). 
\end{align*}
The intrinsic reason is that the general procedure does not assume any structure of the family of candidate estimators. Hence advanced multiscale theory is not applicable. 

\paragraph{Remark.}
In case of unknown $\sigma$, let $\alpha' := 1 - (1 - \alpha)^{1/2}$. Then with probability at least $1 - \alpha'$,
$$
	\alpha'/2 \ \le \ F_m \bigl( m (\hat{\sigma}_n/\sigma)^2 \,\big|\, 0 \bigr)
		\ \le \ 1 - \alpha'/2 .
$$
The latter inequalities entail that $(\sigma/\hat{\sigma}_n)^2$ lies between $\tau_{n,\alpha,l} := m / \chi_{m;1 - \alpha'/2}$ and $\tau_{n,\alpha,u} := m / \chi_{m;\alpha'/2}^2$. Then we obtain simultaneous $(1 - \alpha)$-confidence bounds $\hat{\delta}_{n,\alpha,l}^2(J)$ and $\hat{\delta}_{n,\alpha,u}^2(J)$ as in (\ref{CC.1L}) and (\ref{CC.1U}) by replacing $\alpha$ with $\alpha'$ and $T_n(J)$ with
$$
	\frac{\tau_{n,\alpha,l}}{\hat{\sigma}_n^2} \sum_{i \in J} X_{in}^2
	\quad\text{and}\quad
	\frac{\tau_{n,\alpha,u}}{\hat{\sigma}_n^2} \sum_{i \in J} X_{in}^2 ,
$$
respectively. The conclusions of Theorem~\ref{thm: oracle general} continue to hold, as long as $n/m_n = O(1)$.

%==================
\section{Proofs}
\label{sec: proofs}
%==================

%===================================================
\subsection{Proof of (\ref{eq: performance naive S})
	and (\ref{eq: performance advanced S})}
%===================================================

Note first that $\min_{[0,1]} Y$ lies between $F_n(s_o) + \min_{[0,1]} W$ and $F_n(s_o) + W(s_o)$. Hence for any $\alpha' \in (0,1)$,
\begin{eqnarray*}
	\hat{\SS}_{\alpha}^{\rm naive}
	& \subset & \bigl\{ s \in [0,1] : 
		F_n(s) + W(s) \le F_n(s_o) + W(s_o) + \kappa_\alpha^{\rm naive} \bigr\} \\
	& \subset & \bigl\{ s \in [0,1] : 
		F_n(s) - F_n(s_o) \le \kappa_{\alpha'}^{\rm naive} + \kappa_\alpha^{\rm naive}
		\bigr\} \\
	& = & \bigl\{ s \in [0,1] : 
		F_o(s) - F_o(s_o) \le c_n^{-1} \bigl( \kappa_{\alpha'}^{\rm naive} + \kappa_\alpha^{\rm naive} \bigr)
		\bigr\}
\end{eqnarray*}
and
\begin{eqnarray*}
	\hat{\SS}_{\alpha}^{\rm naive}
	& \supset & \bigl\{ s \in [0,1] : 
		F_n(s) + W(s) \le F_n(s_o) + \min_{[0,1]} W + \kappa_\alpha^{\rm naive} \bigr\} \\
	& \supset & \bigl\{ s \in [0,1] : 
		F_n(s) - F_n(s_o) \le \kappa_\alpha^{\rm naive} - \kappa_{\alpha'}^{\rm naive}
		\bigr\} \\
	& = & \bigl\{ s \in [0,1] : 
		F_o(s) - F_o(s_o) \le c_n^{-1}
			\bigl( \kappa_\alpha^{\rm naive} - \kappa_{\alpha'}^{\rm naive} \bigr)
		\bigr\}
\end{eqnarray*}
with probability $1 - \alpha'$. Since $\kappa_{\alpha'}^{\rm naive} < \kappa_\alpha^{\rm naive}$ if $\alpha < \alpha' < 1$, these considerations, combined with the expansion of $F_o$ near $s_o$, show that the maximum of $|s - s_o|$ over all $s \in \hat{\SS}_\alpha^{\rm naive}$ is precisely of order $O_p(c_n^{-1/\gamma})$.

On the other hand, the confidence region $\hat{\SS}_\alpha$ is contained in the set of all $s \in [0,1]$ such that
$$
		F_n(s) + W(s) \le F_n(s_o) + W(s_o) + \sqrt{|s - s_o|}
			\Bigl( \sqrt{2 \log(e/|s - s_o|)} + \kappa_\alpha \Bigr) \Bigr\} ,
$$
and this entails that
$$
	F_o(s) - F_o(s_o) \le c_n^{-1} \sqrt{|s - s_o|}
			\Bigl( \sqrt{2 \log(e/|s - s_o|)} + \kappa_\alpha + O_p(1) \Bigr)
$$
with $O_p(1)$ not depending on $s$. Now the expansion of $F_o$ near $s_o$ entails claim (\ref{eq: performance advanced S}).	\hfill	$\Box$

%====================================
\subsection{Exponential inequalities}
%====================================

An essential ingredient for our main results is an exponential inequality for quadratic functions of a Gaussian random vector. It extends inequalities of Dahlhaus and Polonik (2006) for quadratic forms and is of independent interest.

\begin{proposition}
\label{prop: exponential chi2}
Let $Z_1, \ldots, Z_n$ be independent, standard Gaussian random variables. Furthermore, let $\lambda_1, \ldots, \lambda_n$ and $\delta_1, \ldots, \delta_n$ be real constants, and define $\gamma^2 := \mathrm{Var} \Bigl( \sum_{i=1}^n \lambda_i(Z_i + \delta_i)^2 \Bigr) = \sum_{i=1}^n \lambda_i^2(2 + 4\delta_i^2)$. Then for arbitrary $\eta \ge 0$ and $\lambda_{\rm max} := \max(\lambda_1,\ldots,\lambda_n,0)$,
\begin{eqnarray*}
	\Pr \Bigl( \sum_{i=1}^n \lambda_i \bigl( (Z_i + \delta_i)^2 - (1 + \delta_i^2) \bigr)
		\ge \eta\gamma \Bigr)
	& \le & \exp \Bigl(- \, \frac{\eta^2}{2 + 4 \eta \lambda_{\rm max}/\gamma} \Bigr) \\
	& \le & e^{1/4} \exp \bigl( - \eta/\sqrt{8} \bigr) .
\end{eqnarray*}
\end{proposition}

Note that replacing $\lambda_i$ in Proposition~\ref{prop: exponential chi2} with $- \lambda_i$ yields twosided exponential inequalities. By means of Proposition~\ref{prop: exponential chi2} and elementary calculations one obtains exponential and related inequalities for noncentral $\chi^2$ distributions:

\begin{corollary}
\label{cor: exponential chi2}
For an integer $n > 0$ and a constant $\delta \ge 0$ let $F_n(\cdot \mid \delta^2)$ be the distribution function of $\chi_n^2(\delta^2)$. Then for arbitrary $r \ge 0$,
\begin{eqnarray}
	F_n(n + \delta^2 + r \mid \delta^2)
	& \ge & 1 - \exp \Bigl( - \, \frac{r^2}{4n + 8\delta^2 + 4r} \Bigr) ,
	\label{cdf noncentral chi2.1}\\
	F_n(n + \delta^2 - r \mid \delta^2)
	& \le & \exp \Bigl( - \, \frac{r^2}{4n + 8\delta^2} \Bigr) .
	\label{cdf noncentral chi2.2}
\end{eqnarray}
In particular, for any $u \in (0,1/2)$,
\begin{eqnarray}
	F_n^{-1}(1 - u \mid \delta^2)
	& \le & n + \delta^2 + \sqrt{(4n + 8\delta^2) \log(u^{-1})} + 4\log(u^{-1}) ,
	\label{quantiles noncentral chi2.1} \\
	F_n^{-1}(u \mid \delta^2)
	& \ge & n + \delta^2 - \sqrt{(4n + 8\delta^2) \log(u^{-1})} .
	\label{quantiles noncentral chi2.2}
\end{eqnarray}
Moreover, for any number $\hat{\delta} \ge 0$, the inequalities $u \le F_n(n + \hat{\delta}^2 \mid \delta^2) \le 1 - u$ entail that
\begin{equation}
\label{confidence bounds noncentral chi2}
	\delta^2 - \hat{\delta}^2 \ \begin{cases}
		\le \ + \sqrt{(4n + 8\hat{\delta}^2) \log(u^{-1})} + 8\log(u^{-1}) , \\
		\ge \ - \sqrt{(4n + 8\hat{\delta}^2) \log(u^{-1})} .
	\end{cases}
\end{equation}
\end{corollary}

Conclusion (\ref{confidence bounds noncentral chi2}) follows from (\ref{cdf noncentral chi2.1}) and (\ref{cdf noncentral chi2.2}), applied to $r = \hat{\delta}^2 - \delta^2$ and $r = \delta^2 - \hat{\delta}^2$, respectively.

\paragraph{Proof of Proposition~\ref{prop: exponential chi2}.}
Standard calculations show that for $0 \le t < (2 \lambda_{\rm max})^{-1}$,
$$
	\mathbb{E} \exp \Bigl( t\sum_{i=1}^n \lambda_i (Z_i + \delta_i)^2 \Bigr)
	\ = \ \exp \Bigl( \frac{1}{2} \sum_{i=1}^n
		\Bigl\{ \delta_i^2 \frac{2t\lambda_i}{1 - 2t\lambda_i} - \log(1 - 2t\lambda_i) \Bigr\} \Bigr) .
$$
Then for any such $t$,
\begin{eqnarray}
	\lefteqn{ \Pr \Bigl( \sum_{i=1}^n \lambda_i \bigl( (Z_i + \delta_i)^2 - (1 + \delta_i^2) \bigr)
		\ge \eta\gamma \Bigr) }
	\nonumber\\
	& \le & \exp \Bigl(- t\eta\gamma - t \sum_{i=1}^n \lambda_i (1 + \delta_i^2) \Bigr)
		\cdot \mathbb{E} \exp \Bigl( t \sum_{i=1}^n \lambda_i (Z_i + \delta_i)^2 \Bigr)
	\nonumber\\
	& = & \exp \Bigl( - t \eta\gamma + \frac{1}{2} \sum_{i=1}^n
		\Bigl\{ \delta_i^2 \frac{4t^2\lambda_i^2}{1 - 2t\lambda_i}
			- \log(1 - 2t\lambda_i) - 2t\lambda_i \Bigr\} \Bigr) .
	\label{hallo}
\end{eqnarray}
Elementary considerations reveal that
$$
	- \log(1 - x) - x \ \le \ \begin{cases}
		x^2/2 & \text{if} \ x \le 0 , \\
		x^2/(2(1 - x)) & \text{if} \ x \ge 0 .
	\end{cases}
$$
Thus (\ref{hallo}) is not greater than
\begin{align*}
	\exp \Bigl( & - t \eta\gamma + \frac{1}{2} \sum_{i=1}^n
		\Bigl\{ \delta_i^2 \frac{4 t^2 \lambda_i^2}{1 - 2t\lambda_i}
			+ \frac{2 t^2 \lambda_i^2}{1 - 2t \max(\lambda_i,0)} \Bigr\} \Bigr) \\
	& \le \ \exp \Bigl( - t \eta \gamma + \frac{\gamma^2 t^2/2 }{1 - 2t\lambda_{\max}} \Bigr) .
\end{align*}
Setting
$$
	t \ := \ \frac{\eta}{\gamma + 2\eta\lambda_{\max}} \ \in \ \left[ 0, (2\lambda_{\rm max})^{-1} \right) ,
$$
the preceding bound becomes
$$
	\Pr \Bigl( \sum_{i=1}^n \lambda_i \bigl( (Z_i + \delta_i)^2 - (1 + \delta_i^2) \bigr)
		\ge \eta\gamma \Bigr)
	\ \le \ \exp \Bigl( - \, \frac{\eta^2}{2 + 4 \eta \lambda_{\max}/\gamma} \Bigr) .
$$
Finally, since $\gamma \ge \lambda_{\rm max} \sqrt{2}$, the second asserted inequality follows from
$$
	\frac{\eta^2}{2 + 4 \eta \lambda_{\rm max}/\gam}
	\ \ge \ \frac{\eta^2}{2 + \sqrt{8} \eta}
	\ = \ \frac{\eta}{\sqrt{8}} - \frac{\eta}{\sqrt{8} + 4 \eta}
	\ \ge \ \frac{\eta}{\sqrt{8}} - \frac{1}{4} .
	\eqno{\Box}
$$

%======================================
\subsection{Proofs of the main results}
%======================================

Throughout this section we assume without loss of generality that $\sigma = 1$. Further let $\SS_n := \{0,1,\ldots,n\}$ and $\TT_n := \big\{ (j,k) : 0\le j < k\le n \big\}$.

\paragraph{Proof of Theorem~\ref{thm: approximation 2}.}

\paragraph{\sl Step I.}
We first analyze $D_n$ in place of $\hat{D}_n$. To collect the necessary ingredients, let the metric $\rho_n$ on $\TT_n$ pointwise be defined by
\bea
	\rho_n \bigl( (j,k), (j',k') \bigr)
	& := & \sqrt{ \tau_n(j,j')^2 + \tau_n(k,k')^2 }.
\eea
We need bounds for the capacity numbers $\mathrm{D}(u, \TT', \rho_n)$ (cf.\ Section~\ref{sec: appendix}) for certain $u > 0$ and $\TT' \subset \TT$. The proof of Theorem~2.1 of D\"umbgen and Spokoiny (2001) entails that
\begin{equation}
\label{ineq: capacity}
	\mathrm{D} \Bigl( u\delta, \bigl\{t\in\TT_n : \tau_n(t) \le \delta\bigr\}, \rho_n \Bigr)
	\ \le \ 12 u^{-4}\delta^{-2}
	\quad\text{for all} \ u, \delta \in (0,1] .
\end{equation}

Note that for fixed $(j,k) \in \TT_n$, $\pm D_n(j,k)$ may be written as
$$
	\sum_{i=1}^n \lambda_i \bigl( (\eps_{in} + \theta_{in})^2 - (1 + \theta_{in}^2) \bigr)
$$
with
$$
	\lambda_i=\lambda_{in}(j,k)
	\ := \ \pm \big(4\Arrowvert\theta_n\Arrowvert^2+2n\big)^{-1/2} I_{(j,k]}(i) ,
$$
so $|\lambda_i| \le \big(4\Arrowvert\theta_n\Arrowvert^2+2n\big)^{-1/2}$. Hence it follows from Proposition~\ref{prop: exponential chi2} that
$$
	\Pr \Bigl( |D_n(t)| \ge \tau_n(t) \eta \Bigr)
	\ \le \ 2 \exp \Biggl( - \, \frac{\eta^2}
			{2 + 4 \eta  \big(4\Arrowvert\theta_n\Arrowvert^2+2n\big)^{1/2}/ \tau_n(t)}
		\Biggr)
$$
for arbitrary $t \in \TT_n$ and $\eta \ge 0$. One may rewrite this exponential inequality as
\begin{equation}
\label{ineq: exponential sharp}
	\Pr \Bigl( |D_n(t)| \ge \tau_n(t) G_n\big(\eta, \tau_n(t)\big) \Bigr)
	\ \le \ 2 \exp(- \eta)
\end{equation}
for arbitrary $t \in \TT_n$ and $\eta \ge 0$, where
$$
	G_n\big(\eta,\delta\big)
	\ := \ \sqrt{2 \eta} + \frac{4\eta}{\big(4\Arrowvert\theta_n\Arrowvert^2+2n\big)^{1/2}\delta} .
$$
The second exponential inequality in Proposition~\ref{prop: exponential chi2} entails that
\begin{equation}
\label{ineq: exponential rough 1}
	\Pr \Bigl( \bigl| D_n(t) \bigr| \ge \tau_n(t) \eta \Bigr)
	\ \le \ 2 e^{1/4} \exp \bigl( - \eta / \sqrt{8} \bigr)
\end{equation}
and
\begin{equation}
\label{ineq: exponential rough 2}
	\Pr \Bigl( \bigl| D_n(s) - D_n(t) \bigr|
		\ge \sqrt{8} \rho_n(s,t) \eta \Bigr)
	\ \le \ 2e^{1/4} \exp(- \eta)
\end{equation}
for arbitrary $s, t \in \TT_n$ and $\eta \ge 0$.

Utilizing (\ref{ineq: capacity}) and (\ref{ineq: exponential rough 2}), it follows from Theorem~7 and the subsequent Remark~3 in D\"umbgen and Walther (2007) that
\begin{equation}
\label{ineq: modulus}
	\lim_{\delta \downarrow 0} \ \sup_{n} \, \Pr \biggl( \sup_{s,t \in \TT_n : \rho_n(s,t) \le \delta} \,
		\frac{|D_n(s) - D_n(t)|}
		     {\rho_n(s,t) \log\bigl(e/\rho_n(s,t)\bigr)}
		> Q \biggr)
	\ = \ 0
\end{equation}
for a suitable constant $Q > 0$. Since $D_n(j,k) = D_n(0,k) - D_n(0,j)$ and $\tau_n(j,k) = \rho_n \bigl( (0,j), (0,k) \bigr)$, this entails the stochastic equicontinuity of $D_n$ with respect to $\rho_n$. 

For $0 \le \delta < \delta' \le 1$ define
$$
	T_n(\delta,\delta')
	\ := \ \sup_{t\in\TT_n : \delta < \tau_n(t) \le \delta'}
		\Biggl( \frac{|D_n(t)|}{\tau_n(t)} - \Gamma_n(t) -  \frac{c\cdot\Gamma_n(t)^2}{\tau_n(t)\big(4\Arrowvert\theta_n\Arrowvert^2+2n\big)^{1/2}} \Biggr)^+
$$
with a constant $c > 0$ to be specified later. Recall that $\Gamma_n(t) := \big(2 \log\big(e/\tau_n(t)^2\big)^{1/2}$. Starting from (\ref{ineq: capacity}), (\ref{ineq: exponential sharp}) and (\ref{ineq: modulus}), Theorem~8 of D\"umbgen and Walther (2007) and its subsequent remark imply that
\begin{eqnarray}
\label{A.1}
	T_n(0,\delta) \ \to_p \ 0
	\quad\text{as} \ n \to \infty \ \text{and} \ \delta \searrow 0 ,
\end{eqnarray}
provided that $c > 2$. On the other hand, (\ref{ineq: capacity}), (\ref{ineq: exponential rough 1}) and (\ref{ineq: modulus}) entail that
\begin{eqnarray}
\label{A.2}
	T_n(\delta,1) \ = \ O_p(1)
	\quad\text{for any fixed} \ \delta > 0 .
\end{eqnarray}

Now we are ready to prove the first assertion about $\hat{D}_n$. Recall that $\hat{D}_n = \hat{\sigma}_n^{-2} (D_n + V_n)$ and
$$
	\frac{V_n(j,k)}{\tau_n(j,k)}\ =\ \frac{2\beta_n(k-j)}{\tau_n(j,k)\big(4\Arrowvert\theta_n\Arrowvert^2+2n\big)^{1/2}\sqrt{n}}\, Z_n
$$
with $Z_n$ being asymptotically standard normal. Since $\tau_n(j,k) \le \sqrt{2(k - j)}/\big(4\Arrowvert\theta_n\Arrowvert^2+2n\big)^{1/2}$,
\begin{equation}
\label{ineq: easy V}
	\frac{|V_n(j,k)|}{\tau_n(j,k)}
	\ \le \ \frac{\sqrt{2(k - j)}}{\sqrt{n}} \, \beta_n |Z_n|
	\ \le \ \frac{\gamma_n(j,k)}{\sqrt{n}} \, \beta_n |Z_n| ,
\end{equation}
so the maximum of $|V_n|/\tau_n$ over $\TT_n$ is bounded by $\sqrt{2}\beta_n |Z_n| = O_p(1)$. Furthermore, since $|\TT_n| \le n^2/2$, one can easily deduce from (\ref{ineq: exponential rough 1}) that the maximum of $|D_n|/\tau_n$ over $\TT_n$ exceeds $\sqrt{32} \log n + \eta$ with probability at most $e^{1/4} \exp \bigl( - \eta / \sqrt{8} \bigr)$. Since $\hat{\sigma}_n = 1 + O_p(n^{-1/2})$, these considerations show that
$$
	\max_{t \in \TT_n} \frac{\bigl| (D_n + V_n)(t)|}{\tau_n(t)}
	\ \le \ \sqrt{32} \, \log n + O_p(1)
$$
and
$$
	\max_{t \in \TT_n}
		\frac{ \bigl| \hat{D}_n(t) - (D_n + V_n)(t) \bigr|}{\tau_n(t)}
	\ = \ O_p(n^{-1/2} \log n) .
$$
This proves our first assertion about $\hat{D}_n/\tau_n$.

\paragraph{\sl Step II.}
Because $\hat{\sigma}_n^2 \to_p 1$, it is sufficient for the proof of the weak approximation  
\begin{align}
d_w\Big(\big(\hat{D}_n(t)\big)_{t\in\TT_n},\big(\Delta_n(t)\big)_{t\in\TT_n}\Big)\ \rightarrow\ 0 \ \ \ \text{as $n\rightarrow\infty$}\label{WA}
\end{align}
to show the result for 
$
	 \hat{\sigma}_n^2 \hat{D}_n
	\, = \, D_n + V_n
$
with the processes $D_n$ and $V_n$ introduced in (\ref{eq: Definition D_n}) and (\ref{eq: Definition V_n}). Here, $d_w$ refers to the dual bounded Lipschitz metric which metrizes the topology of weak convergence. Further details are provided in the appendix.  
Note that $D_n(j,k) = D_n(k) - D_n(j)$ with $D_n(\ell) := D_n(0,\ell)$ and $V_n(j,k) = V_n(k) - V_n(j)$ with $V_n(\ell) := V_n(0,\ell)$. Thus we view these processes $D_n$ and $V_n$ temporarily as processes on $\SS_n$. They are stochastically independent by Assumption (A). Hence, acccording to Lemma~\ref{lem: BLM for sums}, it suffices to show that $D_n$ and $V_n$ are approximated in distribution by
\begin{equation}
\label{eq: approximating processes}
	\bigl( W(\tau_n(k)) \bigr)_{k \in \SS_n}
	\quad\text{and}\quad
	\biggl( \frac{k}{\sqrt{n} \, \sqrt{4\Arrowvert\theta_n\Arrowvert^2+2n}} \, Z \biggr)_{k \in \SS_n} ,
\end{equation}
respectively. The assertion about $V_n$ is an immediate consequence of the fact that $Z_n := \sqrt{m/2}(1 - \hat{\sigma}_n^2) = \beta_n^{-1} \sqrt{n} (1 - \hat{\sigma}_n^2)$ converges in distribution to $Z$ while $0 \le k /\bigl[\sqrt{n} \big(4\Arrowvert\theta_n\Arrowvert^2+2n\big)^{1/2} \bigr] \le 1 / \sqrt{2}$.

It remains to verify the assertion about $D_n$. It follows from the results in step~I that the sequence of processes $D_n$ on $\SS_n$ is stochastically equicontinuous with respect to the metric $\tau_n$ on $\SS_n\times\SS_n$. 
More precisely, 
$$
	\max_{(j,k) \in \TT_n} \frac{|D_n(k) - D_n(j)|}{\tau_n(j,k) \log(e / \tau_n(j,k)^2)} \ = \ O_p(1) ,
$$
and it is well-known that $\bigl( W(\tau_n(0,k)^2) \bigr)_{k \in \SS_n}$ has the same property, even with the factor $\log(e / \tau_n(j,k)^2)^{1/2}$ in place of $\log(e / \tau_n(j,k)^2)$. Moreover, both processes have independent increments. Thus, in view of Theorem~\ref{thm: weak approximation} in Section~\ref{sec: appendix}, it suffices to show that
\begin{equation}
\label{eq: convergence 1d}
	\max_{(j,k) \in \TT_n} \, d_{\rm w} \Bigl( D_n(j,k), W(\tau_n(0,k)) - W(\tau_n(j)) \Bigr)
	\ \to \ 0 .
\end{equation}

To this end we write $D_n(j,k) = D_{n,1}(j,k) + D_{n,2}(j,k) + D_{n,3}(j,k)$ with
\begin{eqnarray*}
	D_{n,1}(j,k) & := & \big(4\Arrowvert\theta_n\Arrowvert^2+2n\big)^{-1/2} \sum_{i=j+1}^k 1\{|\theta_{in}| \le \delta_n\}
		(2\theta_{in}\eps_{in} + \eps_{in}^2 - 1) , \\
	D_{n,2}(j,k) & := & \big(4\Arrowvert\theta_n\Arrowvert^2+2n\big)^{-1/2} \sum_{i=j+1}^k 1\{|\theta_{in}| > \delta_n\}
		2 \theta_{in}\eps_{in} , \\
	D_{n,3}(j,k) & := & \big(4\Arrowvert\theta_n\Arrowvert^2+2n\big)^{-1/2} \sum_{i=j+1}^k 1\{|\theta_{in}| > \delta_n\}
		(\eps_{in}^2 - 1)
\end{eqnarray*}
and arbitrary numbers $\delta_n > 0$ such that $\delta_n \to \infty$ but $\delta_n / \big(4\Arrowvert\theta_n\Arrowvert^2+2n\big)^{1/2} \to 0$.
These three random variables $D_{n,s}(j,k)$ are uncorrelated and have mean zero. The number $a_n := \bigl| \{i : |\theta_{in}| > \delta_n\} \bigr|$ satisfies the inequality $\|\theta_n\|^2 \ge a_n \delta_n^2$, whence
$$
	\Ex \bigl( D_{n,3}(j,k)^2 \bigr)
	\ \le \ \frac{2 a_n}{2n + 4 \|\theta_n\|^2}
	\ \le \ \frac{1}{2 \delta_n^2}
	\ \to \ 0 .
$$
Moreover, $D_{n,1}(j,k)$ and $D_{n,2}(j,k)$ are stochastically independent, where $D_{n,1}(j,k)$ is asymptotically Gaussian by virtue of Lindeberg's CLT, while $D_{n,2}(j,k)$ is exactly Gaussian. These findings entail (\ref{WA}).

\paragraph{\sl Step III.}
For $0 \le \delta < \delta' \le 1$ define
\begin{align*}
	S_n(\delta,\delta')
	\, &:= \, \sup_{\substack{t \in \TT_n :\\ \delta < \tau_n(t) \le \delta'}}
		\Biggl( \frac{\bigl| (D_n + V_n)(t) \bigr|}{\tau_n(t)} - \Gamma_n(t) -  \frac{c\cdot\Gamma_n(t)^2}{\big(4\Arrowvert\theta_n\Arrowvert^2+2n\big)^{1/2}\tau_n(t)}
		\Biggr)^+ , \\
	\Sigma_n(\delta,\delta')\,
	 &:=\,  \sup_{\substack{(j,k)\in\TT_n :\\ \delta < \tau_n(j,k) \le \delta'}}
		\biggl( \frac{\bigl| W(\tau_n(0,k)^2) - W(\tau_n(0,j)^2) \bigr|}{\tau_n(j,k)}
			- \Gamma_n(j,k) \biggr)^+ .
\end{align*}
Since $S_n(0,1) \le T_n(0,1) + \sqrt{2} \beta_n |Z_n|$, it follows from (\ref{A.1}) and (\ref{A.2}) that $S_n(0,1) = O_p(1)$.

As to the approximation in distribution, since $\tau_n(0,n)\big(4\Arrowvert\theta_n\Arrowvert^2+2n\big)^{1/2} \ge \sqrt{2n} \to \infty$,
$$
	\max_{t : \tau_n(t) \ge \delta} \bigg|\frac{\Gamma_n(t)^2}{\big(4\Arrowvert\theta_n\Arrowvert^2+2n\big)^{1/2}\tau_n(t)}\bigg| \ \to \ 0
	\quad\text{while}\quad
	\max_{t : \tau_n(t) \ge \delta} |\Gamma_n(t)| \ = \ O(1)
$$
for any fixed $\delta \in (0,1)$. Consequently it follows from step II that
\begin{equation}
\label{A}
	d_{\rm w} \bigl( S_n(\delta,1), \Sigma_n(\delta,1) \bigr)
	\ \to \ 0
\end{equation}
for any fixed $\delta \in (0,1)$. Thus it suffices to show that
$$
	S_n(0,\delta), \Sigma_n(0,\delta) \ \to_p \ 0
	\quad\text{as} \ n \to \infty \ \text{and} \ \delta \searrow 0 ,
$$
provided that $\|\theta_n\|^2 = O(n)$. For $\Sigma_n(0,\delta)$ this claim follows, for instance, with the same arguments as (\ref{A.1}). Moreover, $S_n(0,\delta)$ is not greater than
$$
	T_n(0,\delta) + \sup_{t \in\TT_n: \tau_n(t) \le \delta} \frac{|V_n(t)|}{\tau_n(t)}
	\ \le \ T_n(0,\delta) + \frac{\big(4\Arrowvert\theta_n\Arrowvert^2+2n\big)^{1/2}}{\sqrt{n}} \, \delta ,
$$
according to (\ref{ineq: easy V}). Thus our claim follows from (\ref{A.1}) and $ \|\theta_n\|^2 = O(n)$.	\hfill	$\Box$

\paragraph{Proof of Proposition~\ref{prop: coupling}.}
The main ingredient is a well-known representation of noncentral $\chi^2$ distributions as Poisson mixtures of central $\chi^2$ distributions. Precisely,
$$
	\chi_k^2(\delta^2) \ = \ \sum_{j=0}^{\infty} e_{}^{-\delta^2/2} \frac{(\delta^2/2)^j}{j!}
		\cdot \chi_{k+2j}^2 ,
$$
as can be proved via Laplace transforms. Now we define `time
points'
$$
	t_{kn} \ := \ \sum_{i=1}^k \th_{in}^2
	\quad\text{and}\quad
	t_{kn}^* \ := \ t_{j(n)n} + k - j(n)
$$
with $j(n)$ any fixed index in $\KK_n(\th_n)$. This construction entails that $t_{kn}^* \ge t_{kn}$ with equality if, and only if, $k \in \KK_n(\th_n)$.

\begin{figure}
\centerline{\includegraphics[width=0.8\textwidth]{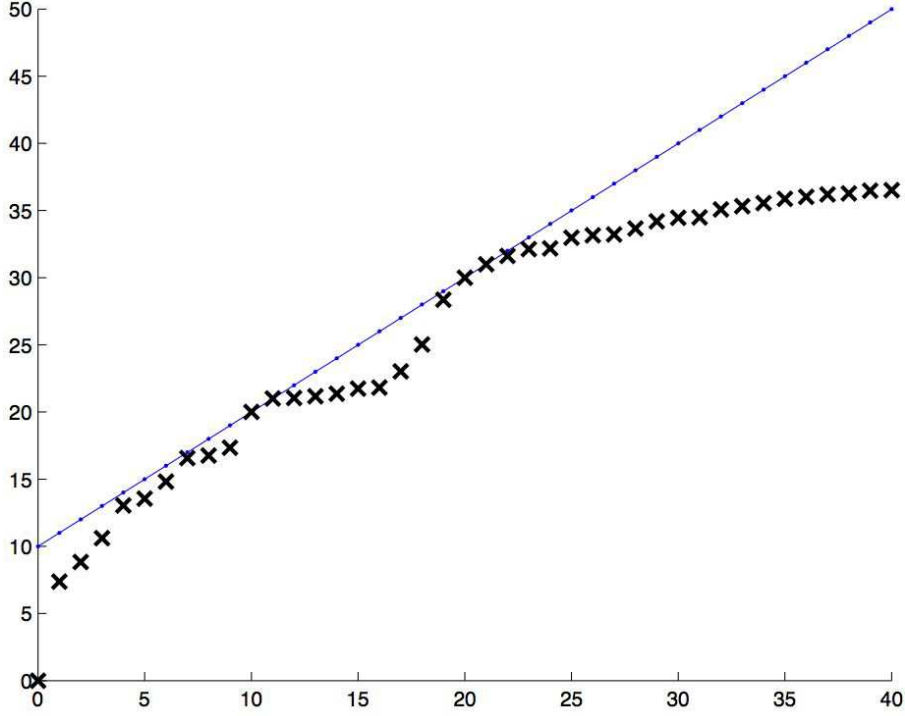}}
\caption{Construction of the coupling}
\label{fig: coupling}
\end{figure}

Figure~\ref{fig: coupling} illustrates this construction. It shows the time points $t_{kn}$ (crosses) and $t_{kn}^*$ (dots and line) versus $k$ for a hypothetical signal $\th_n \in \R^{40}$. Note that in this example, $\KK_n(\th_n)$ is given by $\{10,11, 20, 21\}$.

Let $\Pi$, $G_1$, $G_2$, \ldots, $G_n$, $Z_1$, $Z_2$, $Z_3$, \ldots and $S^2$ be stochastically independent random variables, where $\Pi = (\Pi(t))_{t \ge 0}$ is a standard Poisson process, $G_i$ and $Z_j$ are standard Gaussian random variables, and $S^2 \sim \chi_m^2$. Then one can easily verify that
\bea
	\tilde{T}_{jkn}^{}
	& := & \frac{m}{(k - j) S^2}
		\biggl(\, \sum_{i=j+1}^k G_i^2 + \sum_{s = 2\Pi(t_{jn}/2)+1}^{2 \Pi(t_{kn}/2)} Z_s^2 \biggr) , \\
	\tilde{T}_{jkn}^*
	& := & \frac{m}{(k - j) S^2}
		\biggl(\, \sum_{i=j+1}^k G_i^2 + \sum_{s = 2\Pi(t_{jn}^*/2)+1}^{2 \Pi(t_{kn}^*/2)} Z_s^2 \biggr)
\eea
define random variables $(\tilde{T}_{jkn})_{0 \le j < k \le n}$ and $(\tilde{T}_{jkn}^*)_{0 \le j < k \le n}$ with the desired properties.	\hfill	$\Box$

In the proofs of Theorems~\ref{thm: oracle nested} and \ref{thm: oracle general} we utilize repeatedly two elementary inequalities:

\begin{lemma}
\label{lem: quadratic inequalities}
Let $a,b,c$ be nonnegative constants.

\noindent
\textbf{(i)} \ Suppose that $0 \le x \le y \le x + \sqrt{b(x+y)} + c$. Then
$$
	y \ \le \ x + \sqrt{2bx} + b + \sqrt{bc} + c
	\ \le \ x + \sqrt{2bx} + (3/2)(b + c) .
$$

\noindent
\textbf{(ii)} For $x \ge 0$ define $h(x) := x + \sqrt{a + bx} + c$. Then
$$
	h(h(x)) \ \le \ x + 2 \sqrt{a + bx} + b/2 + \sqrt{bc} + 2c .
$$
\end{lemma}

\paragraph{Proof of Lemma~\ref{lem: quadratic inequalities}.}
The inequality $y \le x + \sqrt{b(x + y)} + c$ entails that either $y < x + c$ or
$$
	(y - x - c)^2 \ \le \ b(x + y) = 2bx + b(y - x) .
$$
Since $y < x + c$ is stronger than the assertions of part~(i), we only consider the displayed quadratic inequality. The latter is equivalent to
$$
	\bigl( y - x - (b/2 + c) \bigr)^2 \ \le \ 2bx + (b/2 + c)^2 - c^2
	\ = \ 2bx + b^2/4 + bc .
$$
Hence the standard inequality $\sqrt{\sum_i z_i} \le \sum_i \sqrt{z_i}$ for nonnegative numbers $z_i$ leads to
$$
	y - x
	\ \le \ \sqrt{2bx} + \sqrt{b^2/4} + \sqrt{bc} + b/2 + c
	\ = \ \sqrt{2bx} + b + \sqrt{bc} + c .
$$
Finally, $0 \le \bigl( \sqrt{b} - \sqrt{c} \bigr)^2$ entails that $\sqrt{bc} \le (b + c)/2$.

As to part~(ii), the definition of $h(x)$ entails that
\bea
	h(h(x))
	& = & x + \sqrt{a + bx}
		+ \sqrt{a + bx + b\sqrt{a + bx} + bc} + 2c \\
	& \le & x + \sqrt{a + bx}
		+ \sqrt{a + bx + b\sqrt{a + bx}} + \sqrt{bc} + 2c \\
	& = & x + \sqrt{a + bx}
		+ \sqrt{a + bx} \sqrt{1 + b/\sqrt{a + bx}} + \sqrt{bc} + 2c \\
	& \le & x + 2 \sqrt{a + bx} + b/2 + \sqrt{bc} + 2c ,
\eea
because $\sqrt{1 + d} \le 1 + d/2$ for arbitrary $d \ge 0$.	\hfill	$\Box$

\paragraph{Proof of Theorem~\ref{thm: oracle nested}.}
The definition of $\hat{\KK}_{n,\alpha}$ and Proposition~\ref{prop: coupling} together entail that $\hat{\KK}_{n,\alpha}$ contains $\KK_n(\theta_n)$ with probability at least $1 - \alpha$. The assertions about $\kappa_{n,\alpha}$ are immediate consequences of Theorem~\ref{thm: approximation 2} applied to $\theta_n = \theta_n^*$.

Now we verify the oracle inequalities (\ref{equ: quality of KKhat1}) and (\ref{equ: quality of KKhat2}). Let $\gamma_n := \bigl(4\|\theta_n\|^2+2n\big)^{1/2} \times \tau_n$. With $\gamma_n^{*}$ we denote the function $\gamma_n$ on $\TT_n$ corresponding to $\theta_n^*$. Throughout this proof we use the shorthand notation $M_n(\ell,k) := M_n(\ell) - M_n(k)$ for $M_n = \hat{R}_n, R_n, \hat{L}_n, L_n$ and arbitrary $\ell,k \in \CC_n$. Furthermore, $\gamma_n^{(*)}(\ell,k) := \gamma_n^{(*)}(k,\ell)$ if $\ell > k$, and $\gamma_n^{(*)}(k,k) := 0$.

In the subsequent arguments, $k_n := \min(\KK_n(\theta_n))$, while $j$ stands for a generic index in $\hat{\KK}_{n,\alpha}$. The definition of the set $\hat{\KK}_{n,\alpha}$ entails that
\begin{equation}
\label{ineq: 1st step for Rn}
	\hat{R}_n(j,k_n)
	\ \le \ \hat{\sigma}_n^2 \Bigl[ \gamma_n^{*}(j,k_n)
		\Bigl( \Gamma \Bigl( \frac{j-k_n}{n} \Bigr) + \kappa_{n,\alpha} \Bigr)
	+ O(\log n) \Bigr] .
\end{equation}
Here and subsequently, $O(r_n)$ and $O_p(r_n)$ denote a generic number and random variable, respectively, depending on $n$ but neither on any other indices in $\CC_n$ nor on $\alpha \in (0,1)$. Precisely, in view of our remark on dependence of $\alpha$, we consider all $\alpha \ge \alpha(n)$ with $\alpha(n) > 0$ such that $\kappa_{n,\alpha(n)} = O(n^{1/6})$. Note that $\hat{\sigma}_n^2 = 1 + O_p(n^{-1/2})$. Moreover, $\gamma_n^{*}(j,k_n)^2 \Gamma \bigl( (j - k_n)/n \bigr)^2$ equals $12 n x \log(e/x) \le 12 n$ with $x := |j - k_n|/n \in [0,1]$. Thus we may rewrite (\ref{ineq: 1st step for Rn}) as
\begin{equation}
\label{ineq: 2nd step for Rn}
	\hat{R}_n(j,k_n)
	\ \le \ \gamma_n^{*}(j,k_n)
		\Bigl( \Gamma \Bigl( \frac{j-k_n}{n} \Bigr) + \kappa_{n,\alpha} \Bigr)
	+ O_p(\log n) .
\end{equation}
Combining this with the equation $R_n(j,k_n) = \hat{R}_n(j,k_n) - \hat{D}_n(j,k_n)$ yields
\begin{equation}
\label{ineq: Rs on confidence region}
	R_n(j,k_n)
	\ \le \ \gamma_n^{*}(j,k_n) \Bigl( \Gamma \Bigl( \frac{j-k_n}{n} \Bigr) + \kappa_{n,\alpha}\Bigr)
		+ O_p(\log n) +  |\hat{D}_n(j,k_n)|.
\end{equation}
Since $\gamma_n^*(j,k_n)^2 \le 6n$ and $\max_{t\in\TT_n} |\hat{D}_n(t)|/\gamma_n(t) = O_p(\log n)$, (\ref{ineq: Rs on confidence region}) yields
$$
	R_n(j,k_n) \ \le \ \sqrt{12 n} + \sqrt{6 n} \, \kappa_{n,\alpha}
		+ O_p(\log n) \gamma_n(j,k_n) .
$$
But elementary calculations yield
\begin{equation}
\label{eq: Rs and gammas}
	\gamma_n(j,k_n)^2
	\ = \ \gamma_n^{*}(j,k_n)^2 + \mathop{\rm sign}(k_n - j) R_n(j,k_n)
	\ \le \ 6 n + R_n(j,k_n) .
\end{equation}
Hence we may conclude that
$$
	R_n(j,k_n)
	\ \le \ O_p(\log n)\sqrt{R_n(j,k_n)}
		+ O_p \bigl( \sqrt{n} (\log n + \kappa_{n,\alpha}) \bigr) ,
$$
and Lemma~\ref{lem: quadratic inequalities}~(i), applied to $x = 0$ and $y = R_n(j,k_n)$, yields
\begin{equation}
\label{ineq: first consistency}
	\max_{j \in \hat{\KK}_{n,\alpha}} R_n(j,k_n)
	\ \le \ O_p \bigl( \sqrt{n} (\log n + \kappa_{n,\alpha}) \bigr) .
\end{equation}

This preliminary result allows us to restrict our attention to indices $j$ in a certain subset of $\CC_n$: Since $0 \le R_n(n,k_n) = n - k_n - \sum_{i=k_n+1}^n \theta_{in}^2$,
$$
	\sum_{i=k_n+1}^n \theta_{in}^2 \ \le \ n - k_n .
$$
On the other hand, in case of $j < k_n$, $R_n(j,k_n) = \sum_{i=j+1}^{k_n} \theta_{in}^2 - (k_n - j)$, so
$$
	\sum_{i=j+1}^n \theta_{in}^2
	\ \le \ n + O_p \bigl( \sqrt{n} (\log n + \kappa_{n,\alpha}) \bigr) .
$$
Thus if $j_n$ denotes the smallest index $j \in \CC_n$ such that $\sum_{i=j+1}^n \theta_{in}^2 \le 2n$, then $k_n \ge j_n$, and $\hat{\KK}_{n,\alpha} \subset \{j_n,\ldots,n\}$ with asymptotic probability one, uniformly in $\alpha \ge \alpha(n)$. This allows us to restrict our attention to indices $j$ in $\{j_n, \ldots, n\} \cap \hat{\KK}_{n,\alpha}$. For any $\ell \ge j_n$, $\hat{D}_n(\ell,k_n)$ involves only the restricted signal vector $(\theta_{in})_{i=j_n+1}^n$, and the proof of Theorem~\ref{thm: approximation 2} entails that
$$
	\max_{j_n \le \ell \le n}
		\biggl( \frac{|\hat{D}_n(\ell,k_n)|}{\gamma_n(\ell,k_n)}
			- \sqrt{2 \log n} - \frac{2c \log n}{\gamma_n(\ell,k_n)} \biggr)^+
	\ = \ O_p(1) .
$$
Thus we may deduce from (\ref{ineq: Rs on confidence region}) the simpler statement that with asymptotic probability one,
\begin{eqnarray}
\label{ineq: Rs on confidence region 2}
	R_n(j,k_n)
	& \le & \bigl( \gamma_n^{*}(j,k_n) + \gamma_n(j,k_n) \bigr)
		\bigl( \sqrt{2 \log n} + \kappa_{n,\alpha} + O_p(1) \bigr) \\
	&& + \ O_p(\log n) .
	\nonumber
\end{eqnarray}
Now we need reasonable bounds for $\gamma_n^{*}(j,k_n)^2$ in terms of $R_n(j)$ and the minimal risk $\rho_n = R_n(k_n)$, where we start from the equation in (\ref{eq: Rs and gammas}): If $j < k_n$, then $\gamma_n(j,k_n)^2 = \gamma_n^{*}(j,k_n)^2 + 4 R_n(j,k_n)$ and $\gamma_n^{*}(j,k_n)^2 = 6 (k_n - j) \le 6 \rho_n$. If $j > k_n$, then $\gamma_n^{*}(j,k_n)^2 = \gamma_n(j,k_n)^2 + 4 R_n(j,k_n)$ and
$$
	\gamma_n(j,k_n)^2 = \sum_{i=k_n+1}^j (4 \theta_{in}^2 + 2)
	\ \le \ 4 \rho_n + 2 R_n(j)
	\ = \ 6 \rho_n + 2 R_n(j,k_n) .
$$
Thus
$$
	\gamma_n^{*}(j,k_n) + \gamma_n(j,k_n)
	\ \le \ 2 \sqrt{6} \sqrt{\rho_n}
		+ \bigl( \sqrt{2} + \sqrt{6} \bigr) \sqrt{R_n(j,k_n)} ,
$$
and inequality (\ref{ineq: Rs on confidence region 2}) leads to
\begin{eqnarray*}
	R_n(j,k_n)
	& \le & \Bigl( 4 \sqrt{3} \sqrt{\log n} + 2 \sqrt{6} \, \kappa_{n,\alpha} + O_p(1) \Bigr)
		\sqrt{\rho_n} \\
	&& + \ O_p \bigl( \sqrt{\log n} + \kappa_{n,\alpha} \bigr) \sqrt{R_n(j,k_n)}
		+ O_p(\log n)
\end{eqnarray*}
for all $j \in \hat{\KK}_{n,\alpha}$. Again we may employ Lemma~\ref{lem: quadratic inequalities} with $x = 0$ and $y = R_n(j,k_n)$ to conclude that
\begin{eqnarray*}
	\max_{j \in \hat{\KK}_{n,\alpha}} \, R_n(j,k_n)
	& \le & \Bigl(4 \sqrt{3} \sqrt{\log n} + 2 \sqrt{6} \, \kappa_{n,\alpha}+O_p(1) \Bigr)
		\sqrt{\rho_n} \\
	&& + \ O_p \Bigl( (\log(n)^{3/4} + \kappa_{n,\alpha(n)}^{3/2}) \rho_n^{1/4}
		+ \log n + \kappa_{n,\alpha(n)}^2 \Bigr)
\end{eqnarray*}
uniformly in $\alpha \ge 0$.

If $\log(n)^3 + \kappa_{n,\alpha(n)}^6 = O(\rho_n)$, then the previous bound for $R_n(j,k_n) = R_n(j) - \rho_n$ reads
$$
	\max_{j \in \hat{\KK}_{n,\alpha}} \, R_n(j)
	\ \le \ \rho_n
		+ \Bigl(4 \sqrt{3} \sqrt{\log n} + 2 \sqrt{6} \, \kappa_{n,\alpha} + O_p(1) \Bigr)
			\sqrt{\rho_n}
$$
uniformly in $\alpha \ge \alpha(n)$. On the other hand, if we consider just a fixed $\alpha > 0$, then $\kappa_{n,\alpha} = O(1)$, and the previous considerations yield
\begin{eqnarray*}
	\max_{j \in \hat{\KK}_{n,\alpha}} \, R_n(j)
	& \le & \rho_n + \bigl(4 \sqrt{3} + o_p(1) \bigr) \sqrt{\log(n) \, \rho_n} \\
	&& + \ O_p \bigl( \log(n)^{3/4} \rho_n^{1/4} + \log n \bigr) \\
	& \le & \rho_n + \bigl(4 \sqrt{3} + o_p(1) \bigr) \sqrt{\log(n) \, \rho_n}
		+ O_p(\log n) .
\end{eqnarray*}
To verify the latter step, note that for any fixed $\epsilon > 0$,
$$
	\log(n)^{3/4} \rho_n^{1/4} \ \le \ \begin{cases}
		\epsilon^{-1} \log n
			& \text{if} \ \rho_n \le \epsilon^{-4} \log n , \\
		\epsilon \sqrt{\log(n) \, \rho_n}
			& \text{if} \ \rho_n \ge \epsilon^{-4} \log n .
	\end{cases}
$$

It remains to prove claim (\ref{equ: quality of KKhat2}) about the losses. From now on, $j$ denotes a generic index in $\CC_n$. Note first that
$$
	L_n(j,k_n) - R_n(j,k_n)
	\ = \ \sum_{i=j+1}^{k_n} (1 - \eps_{in}^2)
	\ = \ R_n(k_n,j) - L_n(k_n,j)
	\quad\text{if} \ j < k .
$$
Thus Theorem \ref{thm: approximation 2}, applied to $\theta_n = 0$, shows that
$$
	\bigl| L_n(j,k_n) - R_n(j,k_n) \bigr|
	\ \le \ \gamma_n^+(j,k_n) \bigl( \sqrt{2 \log n} + O_p(1) \bigr) + O_p(\log n) ,
$$
where
$$
	\gamma_n^+(j,k_n)
	\ := \ \sqrt{2|k_n - j|}
	\ \le \ \sqrt{2 \rho_n} + \sqrt{2 |R_n(j,k)|} .
$$
It follows from $L_n(0) = R_n(0) = \|\theta_n\|^2$ that $L_n(j) - \rho_n$ equals
\bea
	\lefteqn{ L_n(j,k_n) + (L_n - R_n)(k_n,0) } \\
	& = & R_n(j,k_n) + O_p \Bigl( \sqrt{\log(n) \rho_n} \Bigr)
		+ O_p \bigl( \sqrt{\log n} \bigr) \sqrt{R_n(j,k_n)}
		+ O_p(\log n) \\
	& \ge & O_p \Bigl( \sqrt{\log(n) \rho_n} + \log n \Bigr) ,
\eea
because $R_n(j,k_n) \ge 0$ and $R_n(j,k_n) + O_p(r_n) \sqrt{R_n(j,k_n)} \ge O_p(r_n^2)$. Consequently, $\hat{\rho}_n := \min_{j \in \CC_n} L_n(j)$ satisfies the inequality
$$
	\hat{\rho}_n \ \ge \ \rho_n + O_p \Bigl( \sqrt{\log(n) \rho_n} + \log n \Bigr)
	\ = \ (1 + o_p(1)) \rho_n + O_p(\log n) ,
$$
and this is easily shown to entail that
$$
	\rho_n \ \le \ \hat{\rho}_n
		+ O_p \bigl( \sqrt{\log n} \bigr) \sqrt{\hat{\rho}_n} + O_p(\log n)
	\ = \ (1 + o_p(1)) \hat{\rho}_n + O_p(\log n) .
$$

Now we restrict our attention to indices $j \in \hat{\KK}_{n,\alpha}$ again. Here it follows from our result about the maximal risk over $\hat{\KK}_{n,\alpha}$ that $L_n(j) - \rho_n$ equals
\bea
	\lefteqn{ R_n(j,k_n) + O_p \bigl( \sqrt{\log(n) \rho_n} \bigr)
		+ O_p \bigl( \sqrt{\log n} \bigr) \sqrt{R_n(j,k_n)} + O_p(\log n) } \\
	& \le & 2 R_n(j,k_n) + O_p \bigl( \sqrt{\log(n) \rho_n} + \log n \bigr)
		\ \le \ O_p \Bigl( \sqrt{\log(n) \rho_n} + \log n \Bigr) .
\eea
Hence $\max_{j \in \hat{\KK}_{n,\alpha}} L_n(j)$ is not greater than
$$
	\rho_n + O_p \Bigl( \sqrt{\log(n) \rho_n} + \log n \Bigr)
	\ \le \ \hat{\rho}_n + O_p \bigl( \sqrt{\log n} \bigr) \sqrt{\hat{\rho}_n} + O_p(\log n) .
	\eqno{\Box}
$$

\paragraph{Proof of Theorem~\ref{thm: oracle general}.}
The application of inequality (\ref{confidence bounds noncentral chi2}) in Corollary \ref{cor: exponential chi2} to the tripel $(|J|,T_n(J)-|J|, \alpha/(2M_n))$ in place of $(n,\hat{\delta}^2,\alpha)$ yields bounds for $\hat{\delta}_{n,\alpha,l}^2(J)$ and $\hat{\delta}_{n,\alpha,u}^2(J)$ in terms of $\hat{\delta}_n^2(J) := (T_n(J) - |J|)_+$. Then we apply (\ref{quantiles noncentral chi2.1}-\ref{quantiles noncentral chi2.2}) to $T_n(J)$, replacing $(n,\delta^2,u)$ with $(|J|, \delta_n^2(J), \alpha'/(2M_n))$ for any fixed $\alpha' \in (0,1)$. By means of Lemma~\ref{lem: quadratic inequalities}~(ii) we obtain finally
\begin{eqnarray}
\label{Tralala}
	\left.\begin{array}{c}
		\hat{\delta}_{n,\alpha,u}^2(J) - \delta_n^2(J) \\
		\delta_n^2(J) - \hat{\delta}_{n,\alpha,l}^2(J)
	\end{array}\right\}
	& \le & (1+o_p(1))\sqrt{(16 |J| + 32\,\delta_n^2(J)) \log M_n} \\
	&& + \ (K+o_p(1)) \log M_n
	\nonumber	
\end{eqnarray}
for all $J \in \mathcal{M}_n$. Here and throughout this proof, $K$ denotes a generic constant not depending on $n$. Its value may be different in different expressions. It follows from the definition of the confidence region $\hat{\KK}_{n,\alp}$ that for arbitrary $C \in \hat{\KK}_{n,\alpha}$ and $D \in \CC_n$,
\bea
	R_n(C) - R_n(D)
	& = & \delta_n^2(D\setminus C) - \delta_n^2(C\setminus D) + |C| - |D| \\
	& = & (\delta_n^2 - \hat{\delta}_{n,\alp,l}^2)(D\setminus C)
		+ (\hat{\delta}_{n,\alp,u}^2 - \delta_n^2)(C\setminus D) \\
	&& - \ \bigl( \hat{\delta}_{n,\alp,u}^2(C\setminus D)
		- \hat{\delta}_{n,\alp,l}^2(D\setminus C) + |D| - |C| \bigr) \\
	& \le & (\delta_n^2 - \hat{\delta}_{n,\alp,l}^2)(D\setminus C)
		+ (\hat{\delta}_{n,\alp,u}^2 - \delta_n^2)(C\setminus D) .
\eea
Moreover, according to (\ref{Tralala}) the latter bound is not larger than
\bea
	\lefteqn{ (1+o_p(1)) \Bigl\{\sqrt{ \bigl( 16 |D\setminus C| + 32\delta_n^2(D \setminus C) \bigr) \log M_n} } \\
	&& + \ \sqrt{ \bigl( 16 |C\setminus D| + 32\delta_n^2(C \setminus D) \bigr) \log M_n}\Big\}
		+ (K+o_p(1)) \log M_n \\
	& \le & (1+o_p(1)) \sqrt{ 2 \bigl( 16 |D| + 32\delta_n^2(C^{\rm c}) + 16 |C| + 32\delta_n^2(D^{\rm c}) \bigr) \log M_n } \\
	&& \qquad\qquad\qquad\qquad\qquad\qquad + \ (K+o_p(1)) \log M_n \\
	& \le & 8 \sqrt{ \bigl( R_n(C) + R_n(D) \bigr)  \log M_n } \, (1+o_p(1))\ +\ (K+o_p(1)) \log M_n .
\eea
Thus we obtain the quadratic inequality
\bea
	R_n(C) - R_n(D)
	& \le & 8 \sqrt{ \bigl( R_n(C) + R_n(D) \bigr) \log M_n } \, (1+o_p(1)) \\
	&& \qquad\qquad\qquad + \ (K+o_p(1)) \log M_n ,
\eea
and with Lemma~\ref{lem: quadratic inequalities} this leads to
$$
	R_n(C) \ \le \ R_n(D)
		+ 8\sqrt{2} \sqrt{ R_n(D) \log M_n}(1+o_p(1)) + (K + o_p(1)) \log M_n .
$$
This yields the assertion about the risks.

As for the losses, note that $L_n(\cdot)$ and $R_n(\cdot)$ are closely related in that
$$
	(L_n - R_n)(D) \ = \ \sum_{i \in D} \eps_{in}^2 - |J|
$$
for arbitrary $D \in \CC_n$. Hence we may utilize (\ref{quantiles noncentral chi2.1}-\ref{quantiles noncentral chi2.2}), replacing the tripel $(n, \delta^2, u)$ with $(|D|, 0, \alpha'/(2\mu_n))$, to complement (\ref{Tralala}) with the following observation:
\begin{equation}
	- A \sqrt{|D| \log M_n}
	\ \le \ L_n(D) - R_n(D)
	\ \le \ A \sqrt{|D| \log M_n} + A \log M_n
	\label{Tralala2}
\end{equation}
simultaneously for all $D \in \CC_n$ with probability tending to one as $n \to \infty$ and $A \to \infty$. Note also that (\ref{Tralala2}) implies that $R_n(D) \le A \sqrt{R_n(D) \log M_n} + L_n(D)$. Hence
$$
	R_n(D) \ \le \ (3/2) \bigl( L_n(D) + A^2 \log M_n \bigr)
	\quad\text{for all} \ D \in \CC_n ,
$$
by Lemma~\ref{lem: quadratic inequalities}~(i). Assuming that both (\ref{Tralala}) and (\ref{Tralala2}) hold for some large but fixed $A$, we may conclude that for arbitrary $C \in \hat{\KK}_{n,\alpha}$ and $D \in \CC_n$,
\bea
	\lefteqn{ L_n(C) - L_n(D) } \\
	& = & (L_n - R_n)(C) - (L_n - R_n)(D) + R_n(C) - R_n(D) \\
	& \le & A \sqrt{2 (|C| + |D|) \log M_n}
		+ A \sqrt{ 2 \bigl( R_n(C) + R_n(D) \bigr) \log M_n } + 4 A \log M_n \\
	& \le & 2 A \sqrt{ 2 \bigl( R_n(C) + R_n(D) \bigr) \log M_n } + 4 A \log M_n \\
	& \le & A' \sqrt{ \bigl( L_n(C) + L_n(D) \bigr) \log M_n } + A'' \log M_n
\eea
for constants $A'$ and $A''$ depending on $A$. Again this inequality entails that
$$
	L_n(C) \ \le \ L_n(D) + A' \sqrt{2 L_n(D) \log M_n} + A''' \log M_n
$$
for another constant $A''' = A'''(A)$.	\hfill	$\Box$

%==========================
\section{Auxiliary results}
\label{sec: appendix}
%==========================

This section collects some results from the vicinity of empirical process theory which are used in the present paper.

For any pseudo-metric space $(\mathcal{X},d)$ and $u > 0$, we define the capacity number
$$
	\mathrm{D}(u,\mathcal{X},d)
	\ := \ \max \bigl\{ |\mathcal{X}_o| : \mathcal{X}_o \subset \mathcal{X},
		d(x,y) > u \ \text{for different} \ x,y \in \mathcal{X}_o \big\} .
$$

It is well-known that convergence in distribution of random variables with values in a separable metric space may be metrized by the dual bounded Lip\-schitz distance. Now we adapt the latter distance for stochastic processes. Let $\ell_\infty(\TT)$ be the space of bounded functions $x : \TT \to \R$, equipped with supremum norm $\|\cdot\|_\infty$. For two stochastic processes $X$ and $Y$ on $\TT$ with bounded sample paths we define
$$
	d_{\rm w}(X, Y) \ := \ \sup_{f \in \mathcal{H}(\TT)} \, \bigl| \Ex^* f(X) - \Ex^* f(Y) \bigr| ,
$$
where $\Pr^*$ and $\Ex^*$ denote outer probabilities and expectations, respectively, while $\mathcal{H}(\TT)$ is the family of all funtionals $f : \ell_\infty(\TT) \to \R$ such that
$$
	|f(x)| \ \le \ 1	\quad\text{and}\quad
	|f(x) - f(y)| \ \le \ \|x - y\|_\infty	\quad\text{for all} \ x, y \in \ell_\infty(\TT) .
$$

If $d$ is a pseudo-metric on $\TT$, then the modulus of continuity $w(x,\delta \,| \, d)$ of a function $x\in l_{\infty}(\TT)$ is defined as
$$
	w(x,\delta \,|\, d)
	:= \sup_{s,t\in\TT: d(s,t) \le \delta } \, |x(s) - x(t)| .
$$
Furthermore, $\mathcal{C}_u(\TT,d)$ denotes the set of uniformly continuous functions on $(\TT,d)$, that is
$$
	\mathcal{C}_u(\TT,d)
	\ = \ \Big\{ x\in l_{\infty}(\TT): \lim_{\delta\searrow 0} w(x,\delta \,|\, d) = 0 \Big\} .
$$

\begin{theorem}
\label{thm: weak approximation}
For $n = 1, 2, 3, \ldots$ consider stochastic processes $X_n = \bigl(X_n(t)\bigr)_{t\in\TT_n}$ and $Y_n = \bigl(Y_n(t)\bigr)_{t\in\TT_n}$ on a metric space $(\TT_n, \rho_n)$ with bounded sample paths. Then
$$
	d_{\rm w}(X_n, Y_n) \ \to \ 0
$$
provided that the following three conditions are satisfied:

\noindent
(i) \ For arbitrary subsets $\TT_{n,o}$ of $\TT_n$ with $|\TT_{n,o}| = O(1)$,
$$
	d_{\rm w} \Bigl( X_n\big|_{\TT_{n,o}} ,  Y_n\big|_{\TT_{n,o}} \Bigr)
	\ \longrightarrow \ 0 ;
$$
(ii) \ for each number $\eps > 0$,
$$
	\lim_{\delta\searrow 0}
	\underset{n\to\infty}{\lim\sup} \,
		\Pr^* \bigl( w(Z_n, \delta \,|\, \rho_n) > \eps\bigr) \ = \ 0
	\quad\text{for} \ Z_n = X_n, Y_n ;
$$
(iii) \ for any $\delta > 0$,	 \ $\mathrm{D}(\delta,\TT_n,\rho_n) = O(1)$.
\end{theorem}

\paragraph{Proof.}
For any fixed number $\delta > 0$ let $\TT_{n,o}$ be a maximal subset of $\TT_n$ such that $\rho_n(s,t) > \delta$ for differnt $s,t \in \TT_{n,o}$. Then $|\TT_{n,o}| = O(1)$ by Assumption~(iii). Moreover, for any $t \in \TT_n$ there exists a $t_o \in \TT_{n,o}$ such that $\rho_n(t, t_o) \le \delta$. Hence there exists a partition of $\TT_n$ into sets $B_n(t_o)$, $t_o \in \TT_{n,o}$, satisfying $t_o \in B_n(t_o) \subset \bigl\{ t \in \TT_n : \rho_n(t,t_o) \le \delta \bigr\}$. For any function $x$ in $\ell_\infty(\TT_n)$ or $\ell_\infty(\TT_{n,o})$ let $\pi_n x \in \ell_\infty(\TT_n)$ be given by
$$
	\pi_n x(t) \ := \ \sum_{t_o \in \TT_{n,o}} 1\{t \in B_n(t_o)\} x(t_o) .
$$
Then $\pi_n x$ is linear in $x \big|_{\TT_{n,o}}$ with $\|\pi_n x\|_\infty = \bigl\| x \big|_{\TT_{n,o}} \bigr\|_\infty$. Moreover, any $x \in \ell_\infty(\TT_n)$ satisfies the inequality $\| x - \pi_n x \|_\infty \le w(x, \delta \,|\, \rho_n)$. Hence for $Z_n = X_n, Y_n$,
\bea
	d_{\rm w}(Z_n, \pi_nZ_n)
	& \le & \sup_{h \in \HH(\TT_n)} \Ex^* \bigl| h(Z_n) - h(\pi_n Z_n) \bigr| \\
	& \le & \Ex^* \min \bigl( \|Z_n - \pi_n Z_n\|_\infty, 1 \bigr) \\
	& \le & \Ex^* \min \bigl( w(Z_n, \delta \,|\, \rho_n), 1 \bigr) ,
\eea
and this is arbitrarily small for sufficiently small $\delta > 0$ and sufficiently large $n$, according to Assumption~(ii).

Furthermore, elementary considerations reveal that
$$
	d_{\rm w}(\pi_n X_n, \pi_n Y_n)
	\ = \ d_{\rm w} \Bigl( X_n \big|_{\TT_{n,o}} , Y_n \big|_{\TT_{n,o}} \Bigr) ,
$$
and the latter distance converges to zero, because of $|\TT_{n,o}| = O(1)$ and Assumption~(i).

Since
$$
	d_{\rm w}(X_n,Y_n)
	\ \le \ d_{\rm w}(X_n,\pi_nX_n) + d_{\rm w}(Y_n,\pi_nY_n) + d_{\rm w}(\pi_nX_n,\pi_nY_n) ,
$$
these considerations entail the assertion that $d_{\rm w}(X_n,Y_n) \to 0$.	\hfill	$\Box$

Finally, the next lemma provides a useful inequality for $d_{\rm w}(\cdot,\cdot)$ in connection with sums of independent processes.

\begin{lemma}
\label{lem: BLM for sums}
Let $X = X_1 + X_2$ and $Y = Y_1 + Y_2$ with independent random variables $X_1$, $X_2$ and independent random variables $Y_1$, $Y_2$, all taking values in $(\ell_\infty(\TT), \|\cdot\|_\infty)$. Then
$$
	d_{\rm w}(X,Y) \ \le \ d_{\rm w}(X_1,Y_1) + d_{\rm w}(X_2,Y_2) .
$$
\end{lemma}

For this lemma it is important that we consider random variables rather than just stochastic processes with bounded sample paths. Note that a stochastic process on $\TT$ is automatically a random variable with values in $(\ell_\infty(\TT), \|\cdot\|_\infty)$ if (a) the index set $\TT$ is finite, or (b) the process has uniformly continuous sample paths with respect to a pseudo-metric $d$ on $\TT$ such that $N(u,\TT,d) < \infty$ for all $u > 0$.

\paragraph{Proof of Lemma~\ref{lem: BLM for sums}.}
Without loss of generality let the four random variables $X_1$, $X_2$, $Y_1$ and $Y_2$ be defined on a common probability space and stochastically independent. Let $f$ be an arbitrary functional in $\HH(\TT)$. Then it follows from Fubini's theorem that
\begin{eqnarray*}
	\lefteqn{ \bigl| \Ex f(X_1+X_2) - \Ex f(Y_1+Y_2) \bigr| } \\
	& \le & 	\bigl| \Ex f(X_1+X_2) - \Ex f(Y_1+X_2) \bigr|
		+ \bigl| \Ex f(Y_1+X_2) - \Ex f(Y_1+Y_2) \bigr| \\
	& \le & 	\Ex \, \Bigl| \Ex(f(X_1+X_2) \,|\, X_2) - \Ex(f(Y_1+X_2) \,|\, X_2) \Bigr| \\
	&& + \ \Ex \, \Bigl| \Ex(f(Y_1+X_2) \,|\, Y_1) - \Ex(f(Y_1+Y_2) \,|\, Y_1) \Bigr| \\
	& \le & d_{\rm w}(X_1,Y_1) + d_{\rm w}(X_2,Y_2) .
\end{eqnarray*}
The latter inequality follows from the fact that the functionals $x \mapsto f(x + X_2)$ and $x \mapsto f(Y_1 + x)$ belong to $\HH(\TT)$, too. Thus $d_{\rm w}(X,Y) \le d_{\rm w}(X_1,Y_1) + d_{\rm w}(X_2,Y_2)$.	\hfill	$\Box$

\paragraph{Acknowledgement.}
Constructive comments of a referee are gratefully acknowledged.

%============
% References:
%============


\begin{thebibliography}{99}
\bibitem{baraud_04}
	\textsc{Baraud, Y.} (2004).
	\newblock Confidence balls in Gaussian regression.
	\newblock \textsl{Ann.\ Statist.\ \textbf{32}}, 528-551.
\bibitem{beran_96}
	\textsc{Beran, R.} (1996).
	\newblock Confidence sets centered at $C_p$ estimators.
	\newblock \textsl{Ann.\ Inst.\ Statist.\ Math.\ \textbf{48}}, 1-15.
\bibitem{beran_00}
	\textsc{Beran, R.} (2000).
	\newblock REACT scatterplot smoothers: superefficiency through basis economy.
	\newblock \textsl{J.\ Amer.\ Statist.\ Assoc.\ \textbf{95}}, 155-169.
\bibitem{beran_duembgen_98}
	\textsc{Beran, R.} and \textsc{D\"umbgen, L.} (1998).
	\newblock Modulation of estimators and confidence sets.
	\newblock \textsl{Ann.\ Statist.\ \textbf{26}}, 1826-1856.
\bibitem{birge_massart_01}
        \textsc{Birg$\acute{\text{e}}$, L. and Massart, P.} (2001).
        \newblock Gaussian model selection.
        \newblock \textsl{J.\ Eur.\ Math.\ Soc.\ \textbf{3}}, 203-268.
\bibitem{cai_99}
	\textsc{Cai, T.T.} (1999).
	\newblock Adaptive wavelet estimation: a block thresholding and
		oracle inequality approach.
	\newblock \textsl{Ann.\ Statist.\ \textbf{26}}, 1783-1799.
\bibitem{cai_02}
	\textsc{Cai, T.T.} (2002).
	\newblock On block thresholding in wavelet regression: adaptivity,
		block size, and threshold level.
	\newblock \textsl{Statistica Sinica \textbf{12}}, 1241-1273.
\bibitem{cai_low_06}
	\textsc{Cai, T.T.} and \textsc{Low, M.G.} (2006). 
	\newblock Adaptive confidence balls.
	\newblock \textsl{Ann.\ Statist.\ \textbf{34}}, 202-228.
\bibitem{cai_low_07}
	\textsc{Cai, T.T.} and \textsc{Low, M.G.} (2007).
	\newblock Adaptive estimation and confidence intervals for
		convex functions and monotone functions.
	\newblock \textsl{Manuscript in preparation.}
\bibitem{dahlhaus_polonik_06}
	\textsc{Dahlhaus, R.} and \textsc{Polonik, W.} (2006).
	\newblock Nonparametric quasi-maximum likelihood estimation for
		Gaussian locally stationary processes.
	\newblock \textsl{Ann.\ Statist.\ \textbf{34}}, 2790-2824.
\bibitem{donoho_johnstone_94}
	\textsc{Donoho, D.L.} and \textsc{Johnstone, I.M.} (1994).
	\newblock Ideal spatial adaptation by wavelet shrinkage.
	\newblock \textsl{Biometrika\ \textbf{81}}, 425-455.
\bibitem{donoho_johnstone_95}
	\textsc{Donoho, D.L.} and \textsc{Johnstone, I.M.} (1995).
	\newblock Adapting to unknown smoothness via wavelet shrinkage.
	\newblock \textsl{JASA\ \textbf{90}}, 1200-1224.
\bibitem{donoho_johnstone_98}
	\textsc{Donoho, D.L.} and \textsc{Johnstone, I.M.} (1998).
	\newblock Minimax estimation via wavelet shrinkage.
	\newblock \textsl{Ann.\ Statist.\ \textbf{26}}, 879-921.
\bibitem{duembgen_02}
	\textsc{D\"{u}mbgen, L.} (2002).
	\newblock Application of local rank tests to nonparametric regression.
	\newblock \textsl{J.\ Nonpar.\ Statist.\ \textbf{14}}, 511-537.
\bibitem{duembgen_03}
	\textsc{D\"umbgen, L.} (2003).
	\newblock Optimal confidence bands for shape-restricted curves.
	\newblock \textsl{Bernoulli \textbf{9}}, 423-449.
\bibitem{duembgen_spokoiny_01}
	\textsc{D\"{u}mbgen, L.} and \textsc{Spokoiny, V.G.} (2001).
	\newblock Multiscale testing of qualitative hypotheses.
	\newblock \textsl{Ann.\ Statist.\ \textbf{29}}, 124-152.
\bibitem{duembgen_walther_07}
	\textsc{D\"{u}mbgen, L.} and \textsc{Walther, G.} (2007).
	\newblock Multiscale inference about a density.
	\newblock Technical report 56, IMSV, University of Bern.
\bibitem{efromovich_98}
	\textsc{Efromovich, S.} (1998).
	\newblock Simultaneous sharp estimation of functions and their derivatives.
	\newblock \textsl{Ann.\ Statist.\ \textbf{26}}, 273-278.
\bibitem{futschik_99}
	\textsc{Futschik, A.} (1999).
	\newblock Confidence regions for the set of global maximizers of
		nonparametrically estimated curves.
	\newblock \textsl{J.\ Statist.\ Plann.\ Inf.\ \textbf{82}}, 237-250.
\bibitem{genovese_wasserman_05}
	\textsc{Genovese, C.R.} and \textsc{Wassermann, L.} (2005).
	\newblock Confidence sets for nonparametric wavelet regression.
	\newblock \textsl{Ann.\ Statist.\ \textbf{33}}, 698-729.
\bibitem{hengartner_stark_95}
	\textsc{Hengartner, N.W.} and \textsc{Stark, P.B.} (1995). 
	\newblock Finite-sample confidence envelopes for shape-restricted densities.
	\newblock \textsl{Ann.\ Statist.\ \textbf{23}}, 525-550.
\bibitem{hoffmann_lepski_02}
	\textsc{Hoffmann, M.} and \textsc{Lepski, O.} (2002).
	\newblock Random rates in anisotropic regression (with discussion).
	\newblock \textsl{Ann. Statist.\ \textbf{30}}, 325-396.
\bibitem{lepski_etal_97}
	\textsc{Lepski, O.V.}, \textsc{Mammen, E.} and \textsc{Spokoiny, V.G.} (1997).
	\newblock Optimal spatial adaptation to inhomogeneous smoothness:
		an approach based on kernel estimates with variable bandwidth selectors.
	\newblock \textsl{Ann.\ Statist.\ \textbf{25}}, 929-947.
\bibitem{li_89}
	\textsc{Li, K.-C.} (1989).
	\newblock Honest confidence regions for nonparametric regression.
	\newblock \textsl{Ann.\ Statist.\ \textbf{17}}, 1001-1008.
\bibitem{polyak_tsybakov_91}
	\textsc{Polyak, B.T.} and \textsc{Tsybakov, A.B.} (1991).
	\newblock Asymptotic optimality of the $C_p$-test for the
		orthogonal series estimation of regression.
	\newblock \textsl{Theory\ Probab.\ Appl.\ \textbf{35}}, 293-306.
\bibitem{robins_vandervaart_06}
	\textsc{Robins, J.} and \textsc{van der Vaart, A.} (2006).
	\newblock Adaptive nonparametric confidence sets.
	\newblock \textsl{Ann.\ Statist.\ \textbf{34}}, 229-253.
\bibitem{stone_84}
	\textsc{Stone, C.J.} (1984).
	\newblock An asymptotically optimal window selection rule for
		kernel density estimates.
	\newblock \textsl{Ann.\ Statist.\ \textbf{12}}, 1285-1297.
\end{thebibliography}
\end{document}